\documentclass[reqno,11pt]{amsart}
\usepackage{amsfonts}

\usepackage{graphicx}
\usepackage{amsmath}
\usepackage{amsmath}
\usepackage{amsxtra}

\theoremstyle{plain}
\newtheorem{theorem}{Theorem}[section]
\newtheorem{lemma}[theorem]{Lemma}

\newtheorem{corollary}[theorem]{Corollary}

\newtheorem{definition}[theorem]{Definition}

\theoremstyle{definition}

\newtheorem{example}[theorem]{Example}

\numberwithin{equation}{section}
\begin{document}
\title[Generalized Browder's and Weyl's theorems]{Generalized Browder's and
Weyl's theorems \\
for Banach space operators{\normalsize \textbf{\ }}}
\author{Ra\'{u}l E. Curto and Young Min Han}
\thanks{The first named author was partially supported by NSF grants
DMS-0099357 and DMS-0400741. The second named author was supported by Kyung
Hee University Research Fund grant KHU - 20040910}

\begin{abstract}
We find necessary and sufficient conditions for a Banach space operator $T$
to satisfy the generalized Browder's theorem, and we obtain new necessary
and sufficient conditions to guarantee that the spectral mapping theorem
holds for the $B$-Weyl spectrum and for polynomials in $T$. \ We also prove
that the spectral mapping theorem holds for the $B$-Browder spectrum and for
analytic functions on an open neighborhood of $\sigma (T)$. \ As
applications, we show that if $T$ is algebraically $M$-hyponormal, or if $T$
is algebraically paranormal, then the generalized Weyl's theorem holds for $%
f(T)$, where $f\in H((T))$, the space of functions analytic on an open
neighborhood of $\sigma (T)$. \ We also show that if $T$ is reduced by each
of its eigenspaces, then the generalized Browder's theorem holds for $f(T)$,
for each $f\in H(\sigma (T))$.
\end{abstract}

\address{Department of Mathematics, University of Iowa, Iowa City, IA
52242-1419.}
\email{rcurto@math.uiowa.edu}
\address{Department of Mathematics, Kyunghee University, Seoul, Korea
130-701.}
\email{ymhan2004@khu.ac.kr}
\maketitle

\section{\label{sect1}Introduction}

In \cite{Weyl}, H. Weyl proved, for hermitian operators on Hilbert space,
his celebrated theorem on the structure of the spectrum (Equation (\ref{eq11}) below).
\ Weyl's theorem has been extended from hermitian operators to hyponormal
and Toeplitz operators (\cite{Co}), and to several classes of operators
including seminormal operators (\cite{Ber1}, \cite{Ber2}). \ Recently, M.
Berkani and J.J. Koliha \cite{Berkani4} introduced the concepts of
generalized Weyl's theorem and generalized Browder's theorem, and they
showed that $T$ satisfies the generalized Weyl's theorem whenever $T$ is a
normal operator on Hilbert space.

In this paper we extend this result to several classes much larger than that
of normal operators. \ We first find necessary and sufficient conditions for
a Banach space operator $T$ to satisfy the generalized Browder's theorem
(Theorem \ref{thm21}). \ We then characterize the smaller class of operators
satisfying the generalized Weyl's theorem (Theorem \ref{thm24}). \ Next, we
obtain a new necessary and sufficient condition to guarantee that the
spectral mapping theorem holds for the $B$-Weyl spectrum and for polynomials
in $T$ (Theorem \ref{thm210}); this result is then refined in the case when $%
T$ already satisfies the generalized Browder's theorem (Theorem \ref{thm211}%
). \ Along the way we prove that the spectral mapping theorem always holds
for the $B$-Browder spectrum and for analytic functions on an open
neighborhood of $\sigma (T)$ (Theorem \ref{thm29}). \ We have three main
applications of our results: if $T$ is algebraically $M$-hyponormal, or if $T
$ is algebraically paranormal, then the generalized Weyl's theorem holds for 
$f(T)$, for each $f\in H(\sigma (T))$, the space of functions analytic on an
open neighborhood of $\sigma (T)$ (Theorems \ref{thm36} and \ref{algpara},
respectively); and if $T$ is reduced by each of its eigenspaces, then the
generalized Browder's theorem holds for $f(T)$, for each $f\in H(\sigma (T))$
(Corollary \ref{cor219}).

As we shall see below, the concept of Drazin invertibility plays an
important role for the class of $B$-Fredholm operators. \ Let $\mathcal{A}$
be a unital algebra. \ We say that $x\in \mathcal{A}$ is Drazin invertible
of degree $k$ if there exists an element $a\in \mathcal{A}$ such that 
\begin{equation*}
x^{k}ax=x^{k},\;\;\;axa=a,\;\;\;\text{and \ \ \ }xa=ax.
\end{equation*}%
For $a\in \mathcal{A}$, the Drazin spectrum is defined as 
\begin{equation*}
\sigma _{D}(a):=\{\lambda \in \mathbb{C}:a-\lambda \text{ is not Drazin
invertible}\}.
\end{equation*}%
In the case of $T\in \mathcal{B(X)}$, it is well known that $T$ is Drazin
invertible if and only if $T$ has finite ascent and descent, which is also
equivalent to having $T$ decomposed as $T_{1}\oplus T_{2}$, where $T_{1}$ is
invertible and $T_{2}$ is nilpotent. \ 

Throughout this note let $\mathcal{B(X)}$, $\mathcal{B}_{0}\mathcal{(X)}$
and $\mathcal{B}_{00}\mathcal{(X)}$ denote, respectively, the algebra of
bounded linear operators, the ideal of compact operators, and the set of
finite rank operators acting on an infinite dimensional Banach space $%
\mathcal{X}$. \ If $T\in \mathcal{B(X)}$ we shall write $N(T)$ and $R(T)$
for the null space and range of $T$. \ Also, let $\alpha (T):=\dim \;N(T)$, $%
\beta (T):=\dim \;\mathcal{X}/R(T)$, and let $\sigma (T)$, $\sigma _{a}(T)$, 
$\sigma _{p}(T)$, $\sigma _{pi}(T)$, $p_{0}(T)$ and $\pi _{0}(T)$ denote the
spectrum, approximate point spectrum, point spectrum, the eigenvalues of
infinite multiplicity of $T$, the set of poles of $T$, and the set of all
eigenvalues of $T$ which are isolated in $\sigma (T)$, respectively. \ An
operator $T\in \mathcal{B(X)}$ is called \textit{upper semi-Fredholm} if it
has closed range and finite dimensional null space, and is called \textit{%
lower semi-Fredholm} if it has closed range and its range has finite
co-dimension. \ If $T\in \mathcal{B(X)}$ is either upper or lower
semi-Fredholm, then $T$ is called \textit{semi-Fredholm}; the \textit{index}
of a semi Fredholm operator $T\in \mathcal{B(X)}$ is defined as 
\begin{equation*}
i(T):=\alpha (T)-\beta (T).
\end{equation*}%
If both $\alpha (T)$ and $\beta (T)$ are finite, then $T$ is called \textit{%
Fredholm}. $\ T\in \mathcal{B(X)}$ is called \textit{Weyl} if it is Fredholm
of index zero, and \textit{Browder} if it is Fredholm ``of finite ascent and
descent;'' equivalently, (\cite[Theorem 7.9.3]{Har2}) if $T$ is Fredholm and 
$T-\lambda $ is invertible for sufficiently small $\lambda \neq 0$ in $%
\mathbb{C}$. \ The essential spectrum, $\sigma _{e}(T)$, the Weyl spectrum, $%
\omega (T)$, and the Browder spectrum, $\sigma _{b}(T)$, are defined as (%
\cite{Har1},\cite{Har2}) 
\begin{equation*}
\sigma _{e}(T):=\{\lambda \in \mathbb{C}:T-\lambda \ \text{is not Fredholm}%
\},
\end{equation*}%
\begin{equation*}
\omega (T):=\{\lambda \in \mathbb{C}:T-\lambda \ \text{is not Weyl}\},
\end{equation*}%
and 
\begin{equation*}
\sigma _{b}(T):=\{\lambda \in \mathbb{C}:T-\lambda \ \text{is not Browder}\},
\end{equation*}%
respectively. \ Evidently 
\begin{equation*}
\sigma _{e}(T)\subseteq \omega (T)\subseteq \sigma _{b}(T)=\sigma
_{e}(T)\cup \operatorname*{acc}\;\sigma (T),
\end{equation*}%
where we write $\operatorname*{acc}\;K$ for the accumulation points of ${K}%
\subseteq \mathbb{C}$. \ For $T\in \mathcal{B(X)}$ and a nonnegative integer 
$n$ we define $T_{[n]}$ to be the restriction of $T$ to $R(T^{n})$, viewed
as a map from $R(T^{n})$ into $R(T^{n})$ (in particular $T_{[0]}=T$). \ If
for some integer $n$ the range $R(T^{n})$ is closed and $T_{[n]}$ is upper
(resp. lower) semi-Fredholm, then $T$ is called \textit{upper} (resp. 
\textit{lower}) \textit{semi}-$B$-\textit{Fredholm}. \ Moreover, if $T_{[n]}$
is Fredholm, then $T$ is called $B$-Fredholm. \ $T$ is called \textit{semi}-$%
B$-\textit{Fredholm} if it is upper or lower $\text{semi}$-$B$-$\text{%
Fredholm}$.

\begin{definition}
\label{def11}Let $T\in \mathcal{B(X)}$ and let 
\begin{equation*}
\Delta (T):=\{n\in \mathbb{Z}_{+}:m\in \mathbb{Z}_{+},m\geq n\Rightarrow
R(T^{n})\cap N(T\mathbb{)}\subseteq R(T^{m})\cap N(T)\}.
\end{equation*}%
The \textit{degree of stable iteration of }$T$ is defined as $\operatorname*{dis}%
\;T:=\inf \;\Delta (T)$.
\end{definition}

Let $T$ be $\text{semi}$-$B$-$\text{Fredholm}$ and let $d$ be the degree of
stable iteration of $T$. \ It follows from \cite[Proposition 2.1]{Berkani6}
that $T_{[m]}$ is semi-Fredholm and $i(T_{[m]})=i(T_{[d]})$ for every $m\geq
d$. \ This enables us to define the \textit{index} of a \textit{semi}-$B$-%
\textit{Fredholm} operator $T$ as the index of the semi-Fredholm operator $%
T_{[d]}$. \ Let $BF(\mathcal{X})$ be the class of all $B$-Fredholm
operators. \ In \cite{Berkani1} the author studied this class of operators
and proved \cite[Theorem 2.7]{Berkani1} that $T\in \mathcal{B(X)}$ is $B$%
-Fredholm if and only if $T=T_{1}\oplus T_{2}$, where $T_{1}$ is Fredholm
and $T_{2}$ is nilpotent.

An operator $T\in \mathcal{B(X)}$ is called $B$-Weyl if it is $B$-Fredholm
of index $0$. \ The $B$-Fredholm spectrum, $\sigma _{BF}(T)$, and $B$-Weyl
spectrum, $\sigma _{BW}(T)$, are defined as%
\begin{equation*}
\sigma _{BF}(T):=\{\lambda \in \mathbb{C}:T-\lambda \text{ is not }B\text{%
-Fredholm}\}
\end{equation*}%
and 
\begin{equation*}
\sigma _{BW}(T):=\{\lambda \in \mathbb{C}:T-\lambda \text{ is not }B\text{%
-Weyl}\}.
\end{equation*}%
It is well known that the following equality holds \cite{Berkani2}:%
\begin{equation*}
\sigma _{BW}(T)=\bigcap \{\sigma _{D}(T+F):F\in B_{00}(\mathcal{X})\}.
\end{equation*}%
We now introduce the $B$-Browder spectrum $\sigma _{BB}(T)$, defined as 
\begin{equation*}
\sigma _{BB}(T):=\bigcap \{\sigma _{D}(T+F):F\in B_{00}(\mathcal{X})\text{
and }TF=FT\}.
\end{equation*}%
Clearly, $\sigma _{BW}(T)\subseteq \sigma _{BB}(T)$. \ In this note we shall
show that the $B$-Browder spectrum plays an important role in determining
whether an operator satisfies the generalized Browder's theorem.

If we write $\operatorname*{iso}\;K=K\setminus \operatorname*{acc}\;K$ then we let%
\begin{equation*}
\pi _{00}(T):=\{\lambda \in \operatorname*{iso}\;\sigma (T):0<\alpha (T-\lambda
)<\infty \ \}
\end{equation*}%
and 
\begin{equation*}
p_{00}(T):=\sigma (T)\setminus \sigma _{b}(T).
\end{equation*}%
Given $T\in \mathcal{B(X)}$, we say that Weyl's theorem holds for $T$ (or
that $T$ satisfies Weyl's theorem, in symbols, $T\in \mathcal{W}$) if 
\begin{equation}
\sigma (T)\setminus \omega (T)=\pi _{00}(T),  \label{eq11}
\end{equation}%
\text{and that Browder's theorem holds for }$T$\text{ (in symbols, }$T\in 
\mathcal{B}$\text{) if }%
\begin{equation}
\sigma (T)\setminus \omega (T)=p_{00}(T).  \label{eq12}
\end{equation}%
We also say that the generalized Weyl's theorem holds for $T$ (and we write $%
T\in g\mathcal{W}$) if 
\begin{equation}
\sigma (T)\setminus \sigma _{BW}(T)=\pi _{0}(T),  \label{eq13}
\end{equation}%
\text{and that the generalized Browder's theorem holds for }$T$\text{ (in
symbols, }$T\in g\mathcal{B}$\text{) if }%
\begin{equation}
\sigma (T)\setminus \sigma _{BW}(T)=p_{0}(T).  \label{eq14}
\end{equation}%
\text{It is known (\cite{Har3},\cite{Berkani4}) that }%
\begin{equation}
g\mathcal{W}\subseteq g\mathcal{B}\bigcap \mathcal{W}  \label{gw}
\end{equation}%
and that 
\begin{equation}
g\mathcal{B}\bigcup \mathcal{W}\subseteq \mathcal{B}\text{.}  \label{gb}
\end{equation}%
Moreover, given $T\in g\mathcal{B}$, it is clear that $T\in g\mathcal{W}$ if
and only if $p_{0}(T)=\pi _{0}(T)$.

An operator $T\in \mathcal{B(X)}$ is called \textit{isoloid} if every
isolated point of $\sigma (T)$ is an eigenvalue of $T$. \ If $T\in \mathcal{%
B(X)}$, we write $r(T)$ for the spectral radius of $T$; it is well known
that $r(T)\leq ||T||$. \ An operator $T\in \mathcal{B(X)}$ is called \textit{%
normaloid} if $r(T)=||T||$. \ An operator $X\in \mathcal{B(X)}$ is called a
quasiaffinity if it has trivial kernel and dense range. \ An operator $S\in 
\mathcal{B(X)}$ is said to be a quasiaffine transform of $T\in \mathcal{B(X)}
$ (in symbols, $S\prec T$) if there is a quasiaffinity $X\in \mathcal{B(X)}$
such that $XS=TX$. \ If both $S\prec T$ and $T\prec S$, then we say that $S$
and $T$ are quasisimilar.

We say that $T\in \mathcal{B(X)}$ has the \textit{single valued extension
property} (SVEP) at $\lambda _{0}$ if for every open set $U\subseteq \mathbb{%
C}$ containing $\lambda _{0}$ the only analytic solution $f:U\longrightarrow 
\mathcal{X}$ of the equation 
\begin{equation*}
(T-\lambda )f(\lambda )=0\;\;\;(\lambda \in U)
\end{equation*}%
is the zero function (\cite{Fin},\cite{Lau2}). \ An operator $T$ is said to
have SVEP if $T$ has SVEP at every $\lambda \in \mathbb{C}$. \ Given $T\in 
\mathcal{B(X)}$, the \textit{local resolvent set} $\rho _{T}(x)$ of $T$ at
the point $x\in \mathcal{X}$ is defined as the union of all open subsets $%
U\subseteq \mathbb{C}$ for which there is an analytic function $%
f:U\longrightarrow \mathcal{X}$ such that 
\begin{equation*}
(T-\lambda )f(\lambda )=x\ \quad \;(\lambda \in U).
\end{equation*}%
The \textit{local spectrum} $\sigma _{T}(x)$ of $T$ at $x$ is then defined
as 
\begin{equation*}
\sigma _{T}(x):=\mathbb{C\setminus \rho }_{T}(x).
\end{equation*}%
For $T\in \mathcal{B(X)}$, we define the \textit{local}\ (resp. \textit{%
glocal}) \textit{spectral subspaces} of $T$ as follows. \ Given a set $%
F\subseteq \mathbb{C}$ (resp. a closed set $G\subseteq \mathbb{C}$), 
\begin{equation*}
X_{T}(F):=\{x\in \mathcal{X}:\sigma _{T}(x)\subseteq F\}
\end{equation*}%
(resp. 
\begin{align*}
\mathcal{X}_{T}(G)& :=\{x\in \mathcal{X}:\text{there exists an analytic
function } \\
f& :\mathbb{C\setminus }G\rightarrow \mathcal{X}\text{ such that }(T-\lambda
)f(\lambda )=x\text{ for all }\lambda \in \mathbb{C}\setminus G\}).
\end{align*}%
An operator $T\in \mathcal{B(X)}$ has \textit{Dunford's property} (C) if the
local spectral subspace $X_{T}(F)$ is closed for every closed set $%
F\subseteq \mathbb{C}$. \ We also say that $T$ has \textit{Bishop's property}
($\beta $) if for every sequence $f_{n}:U\rightarrow ${$\mathcal{X}$} such
that $(T-\lambda )f_{n}\rightarrow 0$ uniformly on compact subsets in $U$,
it follows that $f_{n}\rightarrow 0$ uniformly on compact subsets in $U$. \
It is well known \cite{Lau1,Lau2} that 
\begin{equation*}
\text{Bishop's property }({\beta })\Longrightarrow \text{Dunford's property }%
(C)\Longrightarrow \text{SVEP}.
\end{equation*}

\section{\label{sect2}Structural Properties of Operators in $g\mathcal{B}$
and $g\mathcal{W}$}

\begin{theorem}
\label{thm21}Let $T\in \mathcal{B(X)}$. \ Then the following statements are
equivalent:\newline
(i) \ $T\in g\mathcal{B}$\text{;}\newline
(ii) \ \text{$\sigma _{BW}(T)=\sigma _{BB}(T)$; }\newline
(iii) \ \text{$\sigma (T)=\sigma _{BW}(T)\cup \pi _{0}(T)$;}\newline
(iv) \ $\operatorname*{acc}\;$\text{$\sigma (T)\subseteq \sigma _{BW}(T)$};\newline
(v) \ \text{$\sigma (T)\setminus \sigma _{BW}(T)\subseteq \pi _{0}(T)$. }
\end{theorem}

\begin{proof}
(i) $\Longrightarrow $ (ii): Suppose that $T\in g\mathcal{B}$. \ Then $%
\sigma (T)\setminus \sigma _{BW}(T)=p_{0}(T)$. \ Let $\lambda \in \sigma
(T)\setminus \sigma _{BW}(T)$; then $\lambda \in p_{0}(T)$, so $T-\lambda $
is Drazin invertible. \ Let $F\in \mathcal{B}_{00}\mathcal{(X)}$ with $TF=FT$%
. \ It follows from \cite[Theorem 2.7]{Berkani3} that $T+F-\lambda $ is also
Drazin invertible. \ Therefore $\lambda \notin \sigma _{D}(T+F)$, and hence $%
\lambda \notin \sigma _{BB}(T)$. \ Thus, $\sigma _{BB}(T)\subseteq \sigma
_{BW}(T)$. \ On the other hand, it follows from \cite[Theorem 4.3]{Berkani2}
that $\sigma _{BW}(T)=\cap \{\sigma _{D}(T+F):F\in \mathcal{B}_{00}\mathcal{%
(X)}\}$. \ Therefore $\sigma _{BW}(T)\subseteq \sigma _{BB}(T)$, and hence $%
\sigma _{BW}(T)=\sigma _{BB}(T)$.

(ii) $\Longrightarrow $ (i): \ We assume that $\sigma _{BW}(T)=\sigma
_{BB}(T)$ and we will establish that $\sigma (T)\setminus \sigma
_{BW}(T)=p_{0}(T)$. \ Suppose first that $\lambda \in \sigma (T)\setminus
\sigma _{BW}(T)$. \ Then $\lambda \in \sigma (T)\setminus \sigma _{BB}(T)$,
and thus there is a finite rank operator $F$ such that $TF=FT$ and $%
T+F-\lambda $ is Drazin invertible, but $T-\lambda $ is not invertible. \
Since $TF=FT$, it follows from \cite[Theorem 2.7]{Berkani3} that $T-\lambda $
is Drazin invertible. \ Therefore $T-\lambda $ has finite ascent and
descent. \ Since $\lambda \in \sigma (T)$, we have $\lambda \in p_{0}(T)$. \
Thus $\sigma (T)\setminus \sigma _{BW}(T)\subseteq p_{0}(T)$.

Conversely, suppose that $\lambda \in p_{0}(T)$. \ Then $T-\lambda $ is
Drazin invertible but not invertible. \ Since $\lambda $ is an isolated
point of $\sigma (T)$, \cite[Theorem 4.2]{Berkani2} implies that $T-\lambda $
is $B$-Weyl. \ Therefore $\lambda \in \sigma (T)\setminus \sigma _{BW}(T)$.
\ Thus $p_{0}(T)\subseteq \sigma (T)\setminus \sigma _{BW}(T)$.

(ii) $\Longrightarrow $ (iii): \ Let $\lambda \in \sigma (T)\setminus \sigma
_{BW}(T)$. \ Then $\lambda \in \sigma (T)\setminus \sigma _{BB}(T)$, and so
there exists a finite rank operator $F$ such that $TF=FT$ and $T+F-\lambda $
is Drazin invertible, but $T-\lambda $ is not invertible. \ Therefore $%
T-\lambda $ is Drazin invertible but not invertible. \ Hence $\lambda \in
\sigma (T)\setminus \sigma _{D}(T)$, and so $\lambda \in \pi _{0}(T)$. \
Thus $\sigma (T)\subseteq \sigma _{BW}(T)\cup \pi _{0}(T)$. \ Since $\sigma
_{BW}(T)\cup \pi _{0}(T)\subseteq \sigma (T)$, always, we must have $\sigma
(T)=\sigma _{BW}(T)\cup \pi _{0}(T)$.

(iii) $\Longrightarrow $ (ii): \ Suppose that $\sigma (T)=\sigma
_{BW}(T)\cup \pi _{0}(T)$. \ To show that $\sigma _{BW}(T)=\sigma _{BB}(T)$
it suffices to show that $\sigma _{BB}(T)\subseteq \sigma _{BW}(T)$. \
Suppose that $\lambda \in \sigma (T)\setminus \sigma _{BW}(T)$. \ Then $%
T-\lambda $ is $B$-Weyl but not invertible. \ Since $\sigma (T)=\sigma
_{BW}(T)\cup \pi _{0}(T)$, we see that $\lambda \in \pi _{0}(T)$. \ In
particular, $\lambda $ is an isolated point of $\sigma (T)$. \ It follows
from \cite[Theorem 4.2]{Berkani2} that $T-\lambda $ is Drazin invertible. \
Therefore $\lambda \notin \sigma _{D}(T)$. \ If $F$ is a finite rank
operator and $FT=TF$ then by \cite[Theorem 2.7]{Berkani3} $\sigma
_{D}(T)=\sigma _{D}(T+F)$. \ Hence $\lambda \notin \sigma _{BB}(T)$, and so $%
\sigma _{BW}(T)=\sigma _{BB}(T)$.

(i) $\Longleftrightarrow $ (iv): \ Suppose that $T\in g\mathcal{B}$. \ Then $%
\sigma (T)\setminus \sigma _{BW}(T)=p_{0}(T)$. \ Let $\lambda \in \sigma
(T)\setminus \sigma _{BW}(T)$. \ Then $\lambda \in p_{0}(T)$, and so $%
\lambda $ is an isolated point of $\sigma (T)$. \ Therefore $\lambda \in
\sigma (T)\setminus \operatorname*{acc}\;\sigma (T)$, and hence $\operatorname*{acc}%
\;\sigma (T)\subseteq \sigma _{BW}(T)$.

Conversely, let $\lambda \in \sigma (T)\setminus \sigma _{BW}(T)$. \ Since $%
\operatorname*{acc}\;\sigma (T)\subseteq \sigma _{BW}(T)$, it follows that $%
\lambda \in \operatorname*{iso}\;\sigma (T)$ and $T-\lambda $ is $B$-Weyl. \ By %
\cite[Theorem 2.3]{Berkani3}, we must have $\lambda \in p_{0}(T)$. \
Therefore $\sigma (T)\setminus \sigma _{BW}(T)\subseteq p_{0}(T)$. \ For the
converse, suppose that $\lambda \in p_{0}(T)$. \ Then $\lambda $ is a pole
of the resolvent of $T$, and so $\lambda $ is an isolated point of $\sigma
(T)$. \ Therefore $\lambda \in \sigma (T)\setminus \operatorname*{acc}\;\sigma (T)$%
. \ It follows from \cite[Theorem 2.3]{Berkani3} that $\lambda \in \sigma
(T)\setminus \sigma _{BW}(T)$. \ Thus $p_{0}(T)\subseteq \sigma (T)\setminus
\sigma _{BW}(T)$, and so $T\in g\mathcal{B}$.

(iv) $\Longleftrightarrow $ (v): Suppose that $\operatorname*{acc}\;\sigma
(T)\subseteq \sigma _{BW}(T)$, and let $\lambda \in \sigma (T)\setminus
\sigma _{BW}(T)$. \ Then $T-\lambda $ is $B$-Weyl but not invertible. \
Since $\operatorname*{acc}\;\sigma (T)\subseteq \sigma _{BW}(T)$, $\lambda $ is an
isolated point of $\sigma (T)$. \ It follows from \cite[Theorem 2.3]%
{Berkani3} that $\lambda $ is a pole of the resolvent of $T$. \ Therefore $%
\lambda \in \pi _{0}(T)$, and hence $\sigma (T)\setminus \sigma
_{BW}(T)\subseteq \pi _{0}(T)$. \ Conversely, suppose that $\sigma
(T)\setminus \sigma _{BW}(T)\subseteq \pi _{0}(T)$ and let $\lambda \in
\sigma (T)\setminus \sigma _{BW}(T)$. \ Then $\lambda \in \pi _{0}(T)$, and
so $\lambda $ is an isolated point of $\sigma (T)$. \ Therefore $\lambda \in
\sigma (T)\setminus \operatorname*{acc}\;\sigma (T)$, which implies that $\operatorname*{%
acc}\;\sigma (T)\subseteq \sigma _{BW}(T)$.
\end{proof}

\begin{corollary}
\label{cor22}Let $T$ be quasinilpotent or algebraic. \ Then $T\in g\mathcal{B%
}$.
\end{corollary}

\begin{proof}
Straightforward from Theorem \ref{thm21} and the fact that $\operatorname*{acc}%
\;\sigma (T)=\emptyset $ whenever $T$ is quasinilpotent or algebraic.
\end{proof}

Recall that $g\mathcal{W}\subseteq g\mathcal{B}$ (cf. (\ref{gw})). \
However, the reverse inclusion does not hold, as the following example shows.

\begin{example}
\label{ex23}Let $\mathcal{X}=\ell _{p}$, let $T_{1},T_{2}\in \mathcal{B(X)}$
be given by 
\begin{equation*}
T_{1}(x_{1},x_{2},x_{3},\cdots ):=(0,\frac{1}{2}x_{1},\frac{1}{3}x_{2},\frac{%
1}{4}x_{3},\cdots )\text{ and }T_{2}:=0,
\end{equation*}%
and let 
\begin{equation*}
T:=%
\begin{pmatrix}
T_{1} & 0 \\ 
0 & T_{2}%
\end{pmatrix}%
\in \mathcal{B(X}\oplus \mathcal{X)}.
\end{equation*}%
Then\ 
\begin{equation*}
\sigma (T)=\omega (T)=\sigma _{BW}(T)=\pi _{0}(T)=\{0\}
\end{equation*}%
and 
\begin{equation*}
p_{0}(T)=\emptyset .
\end{equation*}%
Therefore, $T\in g\mathcal{B}\setminus g\mathcal{W}$.
\end{example}

The next result gives simple necessary and sufficient conditions for an
operator $T\in g\mathcal{B}$ to belong to the smaller class $g\mathcal{W}$.

\begin{theorem}
\label{thm24}Let $T\in g\mathcal{B}$. \ The following statements are
equivalent.\newline
(i) \ $T\in g\mathcal{W}$.\newline
(ii) \ \text{$\sigma _{BW}(T)\cap \pi _{0}(T)=\emptyset $.}$\newline
$(iii) \ \text{$p_{0}(T)=\pi _{0}(T)$.}
\end{theorem}

\begin{proof}
(i) $\Rightarrow $ (ii): Assume $T\in g\mathcal{W}$, that is, $\sigma
(T)\setminus \sigma _{BW}(T)=\pi _{0}(T)$. \ It then follows easily that $%
\sigma _{BW}(T)\cap \pi _{0}(T)=\emptyset $, as required for (ii).

(ii) $\Rightarrow $ (iii): Let $\lambda \in \pi _{0}(T)$. \ The condition in
(ii) implies that $\lambda \in \sigma (T)\setminus \sigma _{BW}(T)$, and
since $T\in g\mathcal{B}$, we must then have $\lambda \in p_{0}(T)$. \ It
follows that $\pi _{0}(T)\subseteq p_{0}(T)$, and since the reverse
inclusion always hold, we obtain (iii).

(iii) $\Rightarrow $ (i): Since $T\in g\mathcal{B}$, we know that $\sigma
(T)\setminus \sigma _{BW}(T)=p_{0}(T)$, and since we are assuming $%
p_{0}(T)=\pi _{0}(T)$, it follows that $\sigma (T)\setminus \sigma
_{BW}(T)=\pi _{0}(T)$, that is, $T\in g\mathcal{W}$.
\end{proof}

Let $T\in \mathcal{B(X)}$ and let $f\in H(\sigma (T))$, where $H(\sigma (T))$
is the space of functions analytic in an open neighborhood of $\sigma (T)$.
\ It is well known that $\omega (f(T))\subseteq f(\omega (T))$ holds. \ The
following theorem shows that a similar result holds for the $B$-Weyl
spectrum. \ To prove this we begin with the following lemma.

\begin{lemma}
\label{lem25}(\cite[Theorem 3.2]{Berkani2}) Let $S$ and $T$ be two commuting 
$B$-Fredholm operators. \ Then the product $ST$ is a $B$-Fredholm operator
and $i(ST)=i(S)+i(T)$.
\end{lemma}

\begin{theorem}
\label{thm26}Let $T\in \mathcal{B(X)}$ and let $f\in H(\sigma (T))$. \ Then 
\begin{equation}
\sigma _{BW}(f(T))\subseteq f(\sigma _{BW}(T)).  \label{261}
\end{equation}
\end{theorem}

\begin{proof}
Observe that if $S$ and $T$ are two commuting $B$-Weyl operators then the
product $ST$ is a $B$-Weyl operator. \ Indeed, suppose that $S$ and $T$ are
both $B$-Weyl. \ Then $S$ and $T$ are both $B$-Fredholm of index $0$. \ It
follows from Lemma \ref{lem25} that $ST$ is $B$-Fredholm and 
\begin{equation*}
i(ST)=i(S)+i(T)=0.
\end{equation*}%
Therefore $ST$ is $B$-Weyl. \ Let $f$ be an analytic function on an open
neighborhood of $\sigma (T)$. Now we show that $\sigma _{BW}(f(T))\subseteq
f(\sigma _{BW}(T))$. \ Suppose that $\lambda \notin f(\sigma _{BW}(T))$. \
Let 
\begin{equation*}
f(T)-\lambda =c_{0}(T-\lambda _{1})(T-\lambda _{2})\cdots (T-\lambda
_{n})g(T),
\end{equation*}%
where $c_{0},\lambda _{1},\lambda _{2},\dots ,\lambda _{n}\in \mathbb{C}$
and $g(T)$ is invertible. \ Since $\lambda \notin f(\sigma _{BW}(T))$, $%
c_{0}(\mu -\lambda _{1})(\mu -\lambda _{2})\cdots (\mu -\lambda _{n})g(\mu
)\neq 0$ for every $\mu \in \sigma _{BW}(T)$. \ Therefore $\mu \neq \lambda
_{i}$ for every $\mu \in \sigma _{BW}(T)$, and hence $T-\lambda _{i}$ is $B$%
-Weyl ($i=1,2,\dots ,n$). \ Since $g(T)$ is invertible, it follows from the
previous observation that 
\begin{equation*}
i(f(T)-\lambda )=\sum_{j=1}^{n}i(T-\lambda _{j})+i(g(T))=0.
\end{equation*}%
Therefore $f(T)-\lambda $ is $B$-Weyl, and hence $\lambda \notin \sigma
_{BW}(f(T))$. \ Thus $\sigma _{BW}(f(T))\subseteq f(\sigma _{BW}(T))$.
\end{proof}

It is well known that $\sigma _{b}(T)=\sigma _{e}(T)\cup \operatorname*{acc}%
\;\sigma (T)$. \ A similar result holds for the $B$-Browder spectrum.

\begin{theorem}
\label{thm27} Let $T\in \mathcal{B(X)}$. \ Then $\sigma _{BB}(T)=\sigma
_{BF}(T)\cup \operatorname*{acc}\;\sigma (T)$.
\end{theorem}

\begin{proof}
Suppose that $\lambda \notin \sigma _{BB}(T)$. \ Since $\sigma _{BB}(T)=\cap
\{\sigma _{D}(T+F):F\in \mathcal{B}_{00}\mathcal{(X)}$ and $TF=FT\}$, there
exists a finite rank operator $F$ such that $TF=FT$ and $\lambda \notin
\sigma _{D}(T+F)$. \ Since $T+F-\lambda $ is Drazin invertible and $TF=FT$,
it follows from \cite[Theorem 2.7]{Berkani3} that $T-\lambda $ is Drazin
invertible. \ Therefore $T-\lambda $ has finite ascent and descent, and
hence $T-\lambda $ can be decomposed as $T-\lambda =T_{1}\oplus T_{2}$,
where $T_{1}$ is invertible and $T_{2}$ is nilpotent. \ It follows from %
\cite[Lemma 4.1]{Berkani2} that $T-\lambda $ is $B$-Fredholm. \ On the other
hand, since $T-\lambda $ has finite ascent and descent, $\lambda $ is an
isolated point of $\sigma (T)$. \ Hence $\lambda \notin \sigma _{BF}(T)\cup 
\operatorname*{acc}\;\sigma (T)$.

Conversely, suppose that $\lambda \notin \sigma _{BF}(T)\cup \operatorname*{acc}%
\;\sigma (T)$. \ Then $T-\lambda $ is $B$-Fredholm and $\lambda $ is an
isolated point of $\sigma (T)$. \ Since $T-\lambda $ is $B$-Fredholm, it
follows from \cite[Theorem 2.7]{Berkani1} that $T-\lambda $ can be
decomposed as $T-\lambda =T_{1}\oplus T_{2}$, where $T_{1}$ is Fredholm and $%
T_{2}$ is nilpotent. $\ $We consider two cases.

\textbf{Case I}.\textbf{\ \ }Suppose that $T_{1}$ is invertible. \ Then $%
T-\lambda $ is Drazin invertible. \ Thus, if $F$ is a finite rank operator
and $TF=FT$, then $T+F-\lambda $ is Drazin invertible by \cite[Theorem 2.7]%
{Berkani3}. \ Therefore $\lambda \notin \sigma _{BB}(T)$.

\textbf{Case II}.\textbf{\ \ }Suppose that $T_{1}$ is not invertible. \ Then 
$0$ is an isolated point of $\sigma (T_{1})$. \ But $T_{1}$ is a Fredholm
operator, hence it follows from the punctured neighborhood theorem that $%
T_{1}$ is Browder. \ Therefore there exists a finite rank operator $S_{1}$
such that $T_{1}+S_{1}$ is invertible and $T_{1}S_{1}=S_{1}T_{1}$. \ Put $%
F:=S_{1}\oplus 0$. \ Then $F$ is a finite rank operator, $TF=FT$ and 
\begin{equation*}
T-\lambda +F=T_{1}\oplus T_{2}+S_{1}\oplus 0=(T_{1}+S_{1})\oplus T_{2}
\end{equation*}%
is Drazin invertible. \ Hence $\lambda \notin \sigma _{BB}(T)$.
\end{proof}

In general, the spectral mapping theorem does not hold for the $B$-Weyl
spectrum, as shown by the following example.

\begin{example}
\label{ex28} Let $U\in B(l_{2})$ be the unilateral shift and consider the
operator 
\begin{equation*}
T:=U\oplus (U^{\ast }+2).
\end{equation*}%
Let $p(z):=z(z-2)$. \ Since $U$ is Fredholm with $i(U)=-1$ and since $U-2$
and $U^{\ast }+2$ are both invertible, it follows that $T$ and $T-2$ are
Fredholm with indices $-1$ and $1$, respectively. \ Therefore $T$ and $T-2$
are both $B$-Fredholm but $T$ is not $B$-Weyl. \ On the other hand, it
follows from the index product theorem that 
\begin{equation*}
i(p(T))=i(T(T-2))=i(T)+i(T-2)=0,
\end{equation*}%
hence $p(T)$ is Weyl. \ Thus $0\notin \sigma _{BW}(p(T))$, whereas $%
0=p(0)\in p(\sigma _{BW}(T))$.
\end{example}

By contrast, the spectral mapping theorem does hold for the Browder spectrum
and analytic functions. \ The following theorem shows that a similar result
holds for the $B$-Browder spectrum.

\begin{theorem}
\label{thm29}Let $T\in \mathcal{B(X)}$ and let $f\in H(\sigma (T))$. \ Then 
\begin{equation*}
\sigma _{BB}(f(T))=f(\sigma _{BB}(T)).
\end{equation*}
\end{theorem}

\begin{proof}
Suppose that $\mu \notin f(\sigma _{BB}(T))$ and set 
\begin{equation*}
h(\lambda ):=f(\lambda )-\mu .
\end{equation*}%
Then $h$ has no zeros in $\sigma _{BB}(T)$. \ Since $\sigma _{BB}(T)=\sigma
_{BF}(T)\cup \operatorname*{acc}\;\sigma (T)$ by Theorem \ref{thm27}, we conclude
that $h$ has finitely many zeros in $\sigma (T)$. \ Now we consider two
cases.

\textbf{Case I}. \ Suppose that $h$ has no zeros in $\sigma (T)$. \ Then $%
h(T)=f(T)-\mu $ is invertible, and so $\mu \notin \sigma _{BB}(f(T))$.

\textbf{Case II}. \ Suppose that $h$ has at least one zero in $\sigma (T)$.
\ Then 
\begin{equation*}
h(\lambda )\equiv c_{0}(\lambda -\lambda _{1})(\lambda -\lambda _{2})\cdots
(\lambda -\lambda _{n})g(\lambda ),
\end{equation*}%
where $c_{0},\lambda _{1},\lambda _{2},\dots ,\lambda _{n}\in \mathbb{C}$
and $g(\lambda )$ is a nonvanishing analytic function on an open
neighborhood. \ Therefore 
\begin{equation*}
h(T)=c_{0}(T-\lambda _{1})(T-\lambda _{2})\cdots (T-\lambda _{n})g(T),
\end{equation*}%
where $g(T)$ is invertible. \ Since $\mu \notin f(\sigma _{BB}(T))$, $%
\lambda _{1},\lambda _{2},\dots ,\lambda _{n}\notin \sigma _{BB}(T)$. \
Therefore $T-\lambda _{i}$ is $B$-Browder, and hence each $T-\lambda _{i}$
is $B$-Weyl ($i=1,2,\dots ,n$). \ But each $\lambda _{i}$ is an isolated
point of $\sigma (T)$, hence it follows from \cite[Theorem 2.3]{Berkani3}
that each $\lambda _{i}$ is a pole of the resolvent of $T$. \ Therefore $%
T-\lambda _{i}$ has finite ascent and descent ($i=1,2,\dots ,n$), so $%
(T-\lambda _{1})(T-\lambda _{2})\cdots (T-\lambda _{n})$ has finite ascent
and descent by \cite[Theorem 7.1]{A.E.Taylor}. \ Since $g(T)$ is invertible, 
$h(T)$ has finite ascent and descent. \ Therefore $h(T)$ is Drazin
invertible, and so $0\notin \sigma _{D}(h(T))$. \ Hence $\mu \notin \sigma
_{BB}(f(T))$. \ It follows from Cases I and II that $\sigma
_{BB}(f(T))\subseteq f(\sigma _{BB}(T))$.

Conversely, suppose that $\lambda \notin \sigma _{BB}(f(T))$. \ Then $%
f(T)-\lambda $ is $B$-Browder. \ We again consider two cases.

\textbf{Case I}. \ Suppose that $f(T)-\lambda $ is invertible. \ Then $%
\lambda \notin \sigma (f(T))=f(\sigma (T))$, and hence $\lambda \notin
f(\sigma _{BB}(T))$.

\textbf{Case II}. \ Suppose that $\lambda \in \sigma (f(T))\setminus \sigma
_{BB}(f(T))$. \ Write 
\begin{equation*}
f(T)-\lambda \equiv c_{0}(T-\lambda _{1})(T-\lambda _{2})\cdots (T-\lambda
_{n})g(T),
\end{equation*}%
where $c_{0},\lambda _{1},\lambda _{2},\dots ,\lambda _{n}\in \mathbb{C}$
and $g(T)$ is invertible. \ Since $f(T)-\lambda $ is $B$-Browder, there is a
finite rank operator $F$ such that $F(f(T)-\lambda )=(f(T)-\lambda )F$ and $%
f(T)-\lambda +F$ is Drazin invertible. \ It follows from \cite[Theorem 2.7]%
{Berkani3} that $f(T)-\lambda $ is Drazin invertible. \ Therefore $%
f(T)-\lambda =c_{0}(T-\lambda _{1})(T-\lambda _{2})\cdots (T-\lambda
_{n})g(T)$ has finite ascent and descent, and hence $T-\lambda _{i}$ has
finite ascent and descent for every $i=1,2,\dots ,n$ by \cite[Theorem 7.1]%
{A.E.Taylor}. \ Therefore each $T-\lambda _{i}$ is Drazin invertible, and so 
$\lambda _{1},\lambda _{2},\dots ,\lambda _{n}\notin \sigma _{BB}(T)$. \ 

We now wish to prove that $\lambda \notin f(\sigma _{BB}(T))$. \ Assume not;
then there exists a $\mu \in \sigma _{BB}(T)$ such that $f(\mu )=\lambda $.
\ Since $g(\mu )\neq 0$, we must have $\mu =\mu _{i}$ for some $i=1,...,n$,
which implies $\mu _{i}\in \sigma _{BB}(T)$, a contradiction. \ Hence $%
\lambda \notin f(\sigma _{BB}(T))$, and so $f(\sigma _{BB}(T))\subseteq
\sigma _{BB}(f(T))$. \ This completes the proof.
\end{proof}

A sufficient condition for the spectral mapping theorem to hold for the $B$%
-Weyl spectrum and analytic functions can be given in terms of the set 
\begin{equation*}
\mathcal{P(X)}:=\{T\in \mathcal{B(X)}:i(T-\lambda )i(T-\mu )\geq 0\ \text{\
for all }\lambda ,\mu \in \mathbb{C}\setminus \sigma _{BF}(T)\}.
\end{equation*}

\begin{theorem}
\label{thm210}Let $T\in \mathcal{B(X)}$. \ Then the following statements are
equivalent:\newline
(i) $\ T\in \mathcal{P(X)}$;\newline
(ii) \ $f(\sigma _{BW}(T))=\sigma _{BW}(f(T))$\ \ for every $f\in H(\sigma
(T))$.
\end{theorem}

\begin{proof}
(i) $\Longrightarrow $ (ii): \ Suppose that $T\in \mathcal{P(X)}$. \ Since $%
\sigma _{BW}(f(T))\subseteq f(\sigma _{BW}(T))$ by Theorem \ref{thm26}, it
suffices to show that $f(\sigma _{BW}(T))\subseteq \sigma _{BW}(f(T))$. \
Suppose that $\lambda \notin \sigma _{BW}(f(T))$ and write 
\begin{equation}
f(T)-\lambda \equiv c_{0}(T-\lambda _{1})(T-\lambda _{2})\cdots (T-\lambda
_{n})g(T),  \label{eq2101}
\end{equation}%
where $c_{0},\lambda _{1},\lambda _{2},\dots ,\lambda _{n}\in \mathbb{C}$
and $g(T)$ is invertible. \ Since $\lambda \notin \sigma _{BW}(f(T))$, the
operator $f(T)-\lambda $ is $B$-Weyl. \ Therefore $f(T)-\lambda $ is $B$%
-Fredholm with index $0$. \ Since the operators on the right-hand side of (%
\ref{eq2101}) commute, it follows from \cite[Corollary 3.3]{Berkani1} that $%
T-\lambda _{i}$ is $B$-Fredholm (all $i=1,2,\dots ,n$). \ Since $T\in 
\mathcal{P(X)}$, $i(T-\lambda )i(T-\mu )\geq 0\ $\ for all $\lambda ,\mu \in 
\mathbb{C}\setminus \sigma _{BF}(T)$. \ So we consider two cases.

\textbf{Case I}. \ Suppose that $i(T-\lambda )\leq 0\ $\ for all $\lambda
\in \mathbb{C}\setminus \sigma _{BF}(T)$. \ Since $f(T)-\lambda $ is $B$%
-Fredholm with index $0$, we have 
\begin{equation}
i(f(T)-\lambda )=\sum_{j=1}^{n}i(T-\lambda _{j})+i(g(T))=0,  \label{eq2102}
\end{equation}%
which implies that $T-\lambda _{i}$ is $B$-Weyl (all $i=1,2,\dots ,n$). \
Therefore $\lambda \notin f(\sigma _{BW}(T))$, and hence $f(\sigma
_{BW}(T))\subseteq \sigma _{BW}(f(T))$.

\textbf{Case II}. \ Suppose that $i(T-\lambda )\geq 0$ for all $\lambda \in 
\mathbb{C}\setminus \sigma _{BF}(T)$. \ Since each $T-\lambda _{i}$ is $B$%
-Fredholm, $i(T-\lambda _{i})\geq 0$ (all $i=1,2,\dots ,n$). \ Since $%
f(T)-\lambda $ is $B$-Fredholm with index $0$, we know that $T-\lambda _{i}$
is $B$-Weyl ($i=1,2,\dots ,n$), by (\ref{eq2102}). \ Therefore $\lambda
\notin f(\sigma _{BW}(T))$, and hence $f(\sigma _{BW}(T))\subseteq \sigma
_{BW}(f(T))$.

From Cases I and II, it follows that $\sigma _{BW}(f(T))=f(\sigma _{BW}(T))$.

(ii) $\Longrightarrow $ (i): \ Suppose that $\sigma _{BW}(f(T))=f(\sigma
_{BW}(T))$ for every $f\in H(\sigma (T))$. \ Assume to the contrary that $%
T\notin \mathcal{P(X)}$. \ Then there exist $\lambda _{1},\lambda _{2}\in 
\mathbb{C}\setminus \sigma _{BF}(T)$ such that $i(T-\lambda _{1})<0$ and $%
i(T-\lambda _{2})>0$. \ Let $m:=-i(T-\lambda _{1})$ and $n:=i(T-\lambda
_{2}) $. \ Define $f(z):=(z-\lambda _{1})^{n}(z-\lambda _{2})^{m}$. $\ $Then 
$f(T)=(T-\lambda _{1})^{n}(T-\lambda _{2})^{m}$ is $B$-Fredholm and 
\begin{equation*}
i(f(T))=i((T-\lambda _{1})^{n}(T-\lambda _{2})^{m})=ni(T-\lambda
_{1})+mi(T-\lambda _{2})=0.
\end{equation*}%
Therefore $f(T)$ is $B$-Weyl, and hence $0\notin \sigma _{BW}(f(T))$. \ On
the other hand, 
\begin{equation*}
0=f(\lambda _{2})\in f(\sigma _{BW}(T))=\sigma _{BW}(f(T)),
\end{equation*}%
a contradiction. \ Hence $T\in \mathcal{P(X)}$.
\end{proof}

In Theorem \ref{thm29}, we proved that the spectral mapping theorem holds
for the $B$-Browder spectrum and analytic functions. \ This might suggest
that the validity of the generalized Browder's theorem for $T$ provides the
right framework for analyzing the equality in (\ref{261}). \ The following
result confirms this.

\begin{theorem}
\label{thm211}Let $T\in \mathcal{B(X)}$. \ Suppose that $T\in g\mathcal{B}$.
\ Then the following statements are equivalent:\newline
(i) \ \text{$T\in \mathcal{P}(\mathcal{X})$; }\newline
(ii) \ \text{$\sigma _{BW}(f(T))=f(\sigma _{BW}(T))$ }for every \text{$f\in
H(\sigma (T))$; }\newline
(iii) \ $f(T)\in g\mathcal{B}$ for every $f\in H(\sigma (T))$.
\end{theorem}

\begin{proof}
(i) $\Longleftrightarrow $ (ii): \ This is straightforward from Theorem \ref%
{thm210}.

(ii) $\Longleftrightarrow $ (iii): \ Suppose that $\sigma
_{BW}(f(T))=f(\sigma _{BW}(T))$ for every $f\in H(\sigma (T))$. \ By Theorem %
\ref{thm29}, 
\begin{equation*}
\sigma _{BB}(f(T))=f(\sigma _{BB}(T))=f(\sigma _{BW}(T))=\sigma _{BW}(f(T)),
\end{equation*}%
whence $f(T)\in g\mathcal{B}$ by Theorem \ref{thm21}.

Conversely, suppose that $f(T)\in g\mathcal{B}$ for every $f\in H(\sigma
(T)) $. \ It follows from Theorem \ref{thm21} that $\sigma
_{BW}(f(T))=\sigma _{BB}(f(T))$. \ By Theorem \ref{thm29}, we have%
\begin{equation*}
\sigma _{BW}(f(T))=\sigma _{BB}(f(T))=f(\sigma _{BB}(T))=f(\sigma _{BW}(T)),
\end{equation*}%
and hence $\sigma _{BW}(f(T))=f(\sigma _{BW}(T))$.
\end{proof}

As a consequence, we obtain the following theorem, which extends a result in %
\cite{Curto1}.

\begin{theorem}
\label{thm212}Let $S,T\in \mathcal{B(X)}$. \ If $T$ has SVEP and $S\prec T$,
then $f(S)\in g\mathcal{B}$ for every $f\in H(\sigma (S))$.
\end{theorem}

\begin{proof}
Suppose that $T$ has SVEP. \ Since $S\prec T$, it follows from the proof of %
\cite[Theorem 3.2]{Curto1} that $S$ has SVEP. \ We now show that $S\in g%
\mathcal{B}$. \ Let $\lambda \in \sigma (S)\setminus \sigma _{BW}(S)$; then $%
S-\lambda $ is $B$-Weyl but not invertible. \ Since $S-\lambda $ is $B$%
-Weyl, it follows from \cite[Lemma 4.1]{Berkani2} that $S-\lambda $ admits
the decomposition $S-\lambda =S_{1}\oplus S_{2},$where $S_{1}$ is Weyl and $%
S_{2}$ is nilpotent. \ Since $S$ has SVEP, $S_{1}$ and $S_{2}$ also have
SVEP. \ Therefore Browder's theorem holds for $S_{1}$, and hence $\omega
(S_{1})=\sigma _{b}(S_{1})$. \ Since $S_{1}$ is Weyl, $S_{1}$ is Browder. \
Hence $\lambda $ is an isolated point of $\sigma (S)$. \ It follows from
Theorem \ref{thm21} that $S\in g\mathcal{B}$. \ 

Now let $f\in H(\sigma (S))$; we shall show that $\sigma
_{BW}(f(S))=f(\sigma _{BW}(S))$. \ To prove this, by Theorem \ref{thm211} it
suffices to show that $i(S-\lambda )i(S-\mu )\geq 0\ $for every $\lambda
,\mu \in \mathbb{C}\setminus \sigma _{BF}(S)$ . \ Let $\lambda ,\mu \in 
\mathbb{C}\setminus \sigma _{BF}(S)$. \ Then $S-\lambda $ and $S-\mu $ are
both $B$-Fredholm, and so it follows from \cite[Theorem 2.7]{Berkani1} that $%
S-\lambda $ and $S-\mu $ can be decomposed as $S-\lambda =S_{1}\oplus S_{2}$
and $S-\mu =S_{3}\oplus S_{4},$where $S_{1}$ and $S_{3}$ are both Fredholm,
and $S_{2}$ and $S_{4}$ are nilpotent. \ Since $S$ has SVEP, $S_{1}$ and $%
S_{3}$ have SVEP. \ By \cite[Theorem 2.6]{Aie}, $S-\lambda $ and $S-\mu $
have finite ascent, which implies $i(S-\lambda )i(S-\mu )\geq 0$. \ Thus $%
\sigma _{BW}(f(S))=f(\sigma _{BW}(S))$. \ It follows from Theorem \ref%
{thm211} that $f(S)\in g\mathcal{B}$.
\end{proof}

We now recall that the generalized Weyl's theorem may not hold for
quasinilpotent operators, and that it does not necessarily transfer to or
from adjoints.

\begin{example}
\label{ex213}On $\mathcal{X}\equiv \ell _{p}$ let%
\begin{equation*}
T(x_{1},x_{2},x_{3},\cdots ):=(\text{$\frac{1}{2}$}x_{2},\text{$\frac{1}{3}$}%
x_{3},\text{$\frac{1}{4}$}x_{4},\cdots ).
\end{equation*}%
Then 
\begin{equation*}
\sigma (T^{\ast })=\sigma _{BW}(T^{\ast })=\{0\}
\end{equation*}%
and 
\begin{equation*}
\pi _{0}(T^{\ast })=\emptyset .
\end{equation*}%
Therefore $T^{\ast }\in g\mathcal{W}$. \ On the other hand, since $\sigma
(T)=\omega (T)=\pi _{00}(T)$, $T\notin \mathcal{W}$. \ Hence $T\notin g%
\mathcal{W}$.
\end{example}

However, the generalized Browder's theorem performs better.

\begin{theorem}
Let $T\in \mathcal{B(X)}$. \ Then the following statements are equivalent:%
\newline
(i) $\ T\in g\mathcal{B}$;\newline
(ii) $\ T^{\ast }\in g\mathcal{B}$.
\end{theorem}

\begin{proof}
Recall that 
\begin{equation*}
\sigma (T)=\sigma (T^{\ast })\text{ and }\sigma _{BW}(T)=\sigma
_{BW}(T^{\ast }).
\end{equation*}%
Therefore, 
\begin{equation*}
\operatorname*{acc}\;\sigma (T)\subseteq \sigma _{BW}(T)\Longleftrightarrow 
\operatorname*{acc}\;\sigma (T^{\ast })\subseteq \sigma _{BW}(T^{\ast }).
\end{equation*}%
It follows from Theorem \ref{thm21} that $T\in g\mathcal{B}$ if and only if $%
T^{\ast }\in g\mathcal{B}$.
\end{proof}

\section{\label{sect3}Operators Reduced by Their Eigenspaces}

Let $\mathcal{H}$ be an infinite dimensional Hilbert space and suppose that $%
T\in \mathcal{B(H)}$ is reduced by each of its eigenspaces. \ If we let 
\begin{equation*}
\mathfrak{M}:=\bigvee \{N(T-\lambda ):\ \lambda \in \sigma _{p}(T)\},
\end{equation*}%
it follows that $\mathfrak{M}$ reduces $T$. \ Let $T_{1}:=T|\mathfrak{M}$
and $T_{2}:=T|\mathfrak{M}^{\perp }$. \ By \cite[Proposition 4.1]{Ber2} we
have:

\begin{itemize}
\item[(i)] $T_{1}$ is a normal operator with pure point spectrum;

\item[(ii)] $\sigma _{p}(T_{1})=\sigma _{p}(T)$;

\item[(iii)] $\sigma (T_{1})=\text{cl}\,\sigma _{p}(T_{1})$ (here cl denotes
closure);

\item[(iv)] $\sigma _{p}(T_{2})=\emptyset $.
\end{itemize}

In \cite[Definition 5.4]{Ber2}, Berberian defined 
\begin{equation*}
\tau (T):=\sigma (T_{2})\cup \operatorname*{acc}\;\sigma _{p}(T)\cup \sigma
_{pi}(T);
\end{equation*}%
we shall call $\tau (T)$ the \textit{Berberian spectrum} of $T$. \ Berberian
proved that $\tau (T)$ is a nonempty compact subset of $\sigma (T)$. \ In
the following theorem we establish a relation amongst the $B$-Weyl, the $B$%
-Browder and the Berberian spectra.

\begin{theorem}
\label{thm215}Suppose that $T\in \mathcal{B(H)}$ is reduced by each of its
eigenspaces. \ Then 
\begin{equation}
\sigma _{BW}(T)=\sigma _{BB}(T)\subseteq \tau (T).  \label{eq2151}
\end{equation}
\end{theorem}

\begin{proof}
Let $\mathfrak{M}$ be the closed linear span of the eigenspaces $N(T-\lambda
)$ ($\lambda \in \sigma _{p}(T)$) and write 
\begin{equation*}
T_{1}:=T|\mathfrak{M}\text{ and }T_{2}:=T|\mathfrak{M}^{\perp }.
\end{equation*}%
From the preceding arguments it follows that $T_{1}$ is normal, $\sigma
_{p}(T_{1})=\sigma _{p}(T)$ and $\sigma _{p}(T_{2})=\emptyset $. \ Toward (%
\ref{eq2151}) we will show that 
\begin{equation}
\sigma _{BW}(T)\subseteq \tau (T)  \label{eq2152}
\end{equation}%
and 
\begin{equation}
\sigma _{BB}(T)\subseteq \sigma _{BW}(T).  \label{eq2153}
\end{equation}%
To establish (\ref{eq2152}) suppose that $\lambda \in \sigma (T)\setminus
\tau (T)$. \ Then $T_{2}-\lambda $ is invertible and $\lambda \in \pi
_{0}(T_{1})$. \ Since $\sigma _{pi}(T)\subseteq \tau (T)$, we see that $%
\lambda \in \pi _{00}(T_{1})$. \ Since $T_{1}$ is normal, it follows from %
\cite[Theorem 4.5]{Berkani2} that $T_{1}\in g\mathcal{W}$. \ Therefore $%
\lambda \in \sigma (T_{1})\setminus \sigma _{BW}(T_{1})$, and hence $%
T-\lambda $ is $B$-Weyl. \ This proves (\ref{eq2152}).

Toward (\ref{eq2153}) suppose that $\lambda \in \sigma (T)\setminus \sigma
_{BW}(T)$. \ Then $T-\lambda $ is $B$-Weyl but not invertible. \ Observe
that if $\mathcal{H}_{1}$ is a Hilbert space and an operator $R\in \mathcal{%
B(H}_{1}\mathcal{)}$ satisfies $\sigma _{BW}(R)=\sigma _{BF}(R)$, then 
\begin{equation}
\sigma _{BW}(R\oplus S)=\sigma _{BW}(R)\cup \sigma _{BW}(S),  \label{eq2154}
\end{equation}%
for every Hilbert space $\mathcal{H}_{2}$ and $S\in \mathcal{B(H}_{2}%
\mathcal{)}$. \ Indeed, if $\lambda \notin \sigma _{BW}(R)\cup \sigma
_{BW}(S)$, then $R-\lambda $ and $S-\lambda $ are both $B$-Weyl. \ Therefore 
$R-\lambda $ and $S-\lambda $ are $B$-Fredholm with index $0$. \ Hence $%
R-\lambda \oplus S-\lambda $ is $B$-Fredholm; moreover, 
\begin{equation*}
i%
\begin{pmatrix}
R-\lambda & 0 \\ 
0 & S-\lambda%
\end{pmatrix}%
=i(R-\lambda )+i(S-\lambda )=0.
\end{equation*}%
Therefore $R\oplus S-\lambda $ is $B$-Weyl, and so $\lambda \notin \sigma
_{BW}(R\oplus S)$, which implies $\sigma _{BW}(R\oplus S)\subseteq \sigma
_{BW}(R)\cup \sigma _{BW}(S)$. \ Conversely, suppose that $\lambda \notin
\sigma _{BW}(R\oplus S)$. \ Then $R\oplus S-\lambda $ is $B$-Fredholm with
index $0$. \ Since $i(R\oplus S-\lambda )=i(R-\lambda )+i(S-\lambda )$ and $%
i(R-\lambda )=0$, we must have $i(S-\lambda )=0$. \ Therefore $R-\lambda $
and $S-\lambda $ are both $B$-Weyl. \ Hence $\lambda \notin \sigma
_{BW}(R)\cup \sigma _{BW}(S)$, which implies $\sigma _{BW}(R)\cup \sigma
_{BW}(S)\subseteq \sigma _{BW}(R\oplus S)$. \ Since $T_{1}$ is normal, we
can now apply (\ref{eq2154}) to $T_{1}$ in place of $R$ to show that $%
T_{1}-\lambda $ and $T_{2}-\lambda $ are both $B$-Weyl. \ But since $\sigma
_{p}(T_{2})=\emptyset $, we see that $T_{2}-\lambda $ is Weyl and injective.
\ Therefore $T_{2}-\lambda $ is invertible, and so $\lambda \in \sigma
(T_{1})\setminus \sigma _{BW}(T_{1})$. \ Since $T_{1}$ is normal, it follows
from \cite[Theorem 4.5]{Berkani2} that $T_{1}\in g\mathcal{W}$, which
implies $\lambda \in \pi _{0}(T_{1})$. \ Hence $\lambda $ is an isolated
point of $\sigma (T_{1})$ and $T_{2}-\lambda $ is invertible. \ Now observe
that if $\mathcal{H}_{1}$ and $\mathcal{H}_{2}$ are Hilbert spaces then the
following equality holds with no other restriction on either $R$ or $S$: 
\begin{equation}
\sigma _{BB}(R\oplus S)=\sigma _{BB}(R)\cup \sigma _{BB}(S),  \label{eq2155}
\end{equation}%
for every $R\in B(\mathcal{H}_{1})$ and $S\in B(\mathcal{H}_{2})$. \ Indeed,
if $\lambda \notin \sigma _{BB}(R)\cup \sigma _{BB}(S)$, then $R-\lambda $
and $S-\lambda $ are both $B$-Browder. \ So $\lambda \notin \sigma _{BB}(R)$
and $\lambda \notin \sigma _{BB}(S)$, and hence there are finite rank
operators $F_{1}$ and $F_{2}$ such that $RF_{1}=F_{1}R$, $SF_{2}=F_{2}S$, $%
R+F_{1}-\lambda $ and $S+F_{2}-\lambda $ are both Drazin invertible. \ Set 
\begin{equation*}
F:=%
\begin{pmatrix}
F_{1} & 0 \\ 
0 & F_{2}%
\end{pmatrix}%
\text{ and }V:=%
\begin{pmatrix}
R & 0 \\ 
0 & S%
\end{pmatrix}%
.
\end{equation*}%
Then $F$ is a finite rank operator such that 
\begin{equation*}
VF=FV\text{ and }V+F-\lambda \equiv 
\begin{pmatrix}
R+F_{1}-\lambda & 0 \\ 
0 & S+F_{2}-\lambda%
\end{pmatrix}%
\end{equation*}%
is Drazin invertible. \ Therefore $\lambda \notin \sigma _{BB}%
\begin{pmatrix}
R & 0 \\ 
0 & S%
\end{pmatrix}%
$, and hence $\sigma _{BB}%
\begin{pmatrix}
R & 0 \\ 
0 & S%
\end{pmatrix}%
\subseteq \sigma _{BB}(R)\cup \sigma _{BB}(S)$. \ Conversely, suppose that $%
\lambda \notin \sigma _{BB}%
\begin{pmatrix}
R & 0 \\ 
0 & S%
\end{pmatrix}%
$. \ It follows from Theorem \ref{thm27} that $%
\begin{pmatrix}
R-\lambda & 0 \\ 
0 & S-\lambda%
\end{pmatrix}%
$ is $B$-Fredholm and $\lambda $ is an isolated point of $%
\begin{pmatrix}
R & 0 \\ 
0 & S%
\end{pmatrix}%
$. \ Since $\sigma 
\begin{pmatrix}
R & 0 \\ 
0 & S%
\end{pmatrix}%
=\sigma (R)\cup \sigma (S)$, it follows that $R-\lambda $ and $S-\lambda $
are both $B$-Fredholm, and $\lambda $ is an isolated point of $\sigma (R)$
and $\sigma (S)$, respectively. \ It follows from Theorem \ref{thm27} that $%
R-\lambda $ and $S-\lambda $ are both $B$-Browder. \ Therefore $\lambda
\notin \sigma _{BB}(R)\cup \sigma _{BB}(S)$, and hence $\sigma _{BB}(R)\cup
\sigma _{BB}(S)\subseteq \sigma _{BB}%
\begin{pmatrix}
R & 0 \\ 
0 & S%
\end{pmatrix}%
$. \ This proves (\ref{eq2155}).

By Theorem \ref{thm27} and (\ref{eq2155}), we have $\lambda \notin \sigma
_{BB}(T)$. \ This proves (\ref{eq2153}) and completes the proof.
\end{proof}

In \cite{Oberai2}, Oberai showed that if $T\in \mathcal{B(X)}$ is isoloid
and if $T\in \mathcal{W}$ then for any polynomial $p$, $p(T)\in \mathcal{W}$
if and only if $\omega (p(T))=p(\omega (T))$. \ We now show that a similar
result holds for the generalized Weyl's theorem. \ We begin with the
following two lemmas, essentially due to Oberai\ \cite{Oberai2}; we include
proofs for the reader's convenience.

\begin{lemma}
\label{lem216}Let $T\in \mathcal{B(X)}$ and let $f\in H(\sigma (T))$. \ Then 
\begin{equation*}
\sigma (f(T))\setminus \pi _{0}(f(T))\subseteq f(\sigma (T)\setminus \pi
_{0}(T)).
\end{equation*}

\begin{proof}
Suppose that $\lambda \in \sigma (f(T))\setminus \pi _{0}(f(T))$. \ By the
spectral mapping theorem, it follows that $\lambda \in f(\sigma
(T))\setminus \pi _{0}(f(T))$. \ We consider two cases.

\textbf{Case I. \ }Suppose that $\lambda $ is not an isolated point of $%
f(\sigma (T))$. \ Then there exists a sequence $\{\lambda _{n}\}\subseteq
f(\sigma (T))$ such that $\lambda _{n}\rightarrow \lambda $. \ Since $%
\lambda _{n}\in f(\sigma (T))$, $\lambda _{n}=f(\mu _{n})$ for some $\mu
_{n}\in \sigma (T)$. \ By the compactness of $\sigma (T)$, there is a
convergent subsequence $\{\mu _{n_{k}}\}$ such that $\mu _{n_{k}}\rightarrow
\mu \in \sigma (T)$. \ It follows that $f(\mu _{n_{k}})\rightarrow \lambda $%
, and therefore $\lambda =f(\mu )$. \ But $\mu \in \sigma (T)\setminus \pi
_{0}(T)$, whence $\lambda =f(\mu )\in f(\sigma (T)\setminus \pi _{0}(T))$.

\textbf{Case II. \ }Suppose now that $\lambda $ is an isolated point of $%
f(\sigma (T))$. \ Since $\lambda \in \pi _{0}(f(T))$ by assumption, it
follows that $\lambda $ cannot be an eigenvalue of $f(T)$. \ Let 
\begin{equation}
f(T)-\lambda =c_{0}(T-\lambda _{1})(T-\lambda _{2})\cdots (T-\lambda
_{n})g(T),  \label{eq2161}
\end{equation}%
where $c_{0},\lambda _{1},\dots ,\lambda _{n}\in \mathbb{C}$ and $g(T)$ is
invertible. \ Since $f(T)-\lambda $ is injective, and the operators on the
right-hand side of (\ref{eq2161}) commute, none of $\lambda _{1},\lambda
_{2},\ldots ,\lambda _{n}$ can be an eigenvalue of $T$. \ Therefore $\lambda
\in f(\sigma (T)\setminus \pi _{0}(T))$.

From Cases I and II we obtain the desired conclusion.
\end{proof}

\begin{lemma}
\label{lem217}Let $T\in \mathcal{B(X)}$ and assume that $T$ is isoloid. \
Then for any $f\in H(\sigma (T))$ we have 
\begin{equation*}
\sigma (f(T))\setminus \pi _{0}(f(T))=f(\sigma (T)\setminus \pi _{0}(T)).
\end{equation*}
\end{lemma}

\begin{proof}
In view of Lemma \ref{lem216} it suffices to prove that $f(\sigma
(T)\setminus \pi _{0}(T))\subseteq \sigma (f(T))\setminus \pi _{0}(f(T))$. \
Suppose that $\lambda \in f(\sigma (T)\setminus \pi _{0}(T))$. \ Then by the
spectral mapping theorem, we must have $\lambda \in \sigma (f(T))$. \ Assume
that $\lambda \in \pi _{0}(f(T))$. \ Then clearly, $\lambda $ is an isolated
point of $\sigma (f(T))$. \ Let 
\begin{equation*}
f(T)-\lambda =c_{0}(T-\lambda _{1})(T-\lambda _{2})\cdots (T-\lambda
_{n})g(T),
\end{equation*}%
where $c_{0},\lambda _{1},\dots ,\lambda _{n}\in \mathbb{C}$ and $g(T)$ is
invertible. \ If for some $i=1,...,n$, $\lambda _{i}\in \sigma (T)$, then $%
\lambda _{i}$ would be an isolated point of $\sigma (T)$. \ But $T$ is
isoloid, hence $\lambda _{i}$ would also be an eigenvalue of $T$. \ Since $%
\lambda \in \pi _{0}(f(T))$, such $\lambda _{i}$ would belong to $\pi
_{0}(T) $. \ Thus, $\lambda =f(\lambda _{i})$ for some $\lambda _{i}\in \pi
_{0}(T)$, and hence $\lambda \in f(\pi _{0}(T))$, a contradiction. \
Therefore $\lambda \notin \pi _{0}(f(T))$, so that $\lambda \in \sigma
(f(T))\setminus \pi _{0}(f(T))$.
\end{proof}

\begin{theorem}
\label{thm218} Suppose that $T\in \mathcal{B(X)}$ is isoloid and $T\in g%
\mathcal{W}$. \ Then for any $f\in H(\sigma (T))$, 
\begin{equation*}
f(T)\in g\mathcal{W}\Longleftrightarrow f(\sigma _{BW}(T))=\sigma
_{BW}(f(T)).
\end{equation*}
\end{theorem}

\begin{proof}
$(\Longrightarrow )$ Suppose $f(T)\in g\mathcal{W}$. \ Then $\sigma
_{BW}(f(T))=\sigma (f(T))\setminus \pi _{0}(f(T))$. \ Since $T$ is isoloid,
it follows from Lemma \ref{lem217} that $f(\sigma (T)\setminus \pi
_{0}(T))=\sigma (f(T))\setminus \pi _{0}(f(T))$. \ But $T\in g\mathcal{W}$,
hence $\sigma _{BW}(T)=\sigma (T)\setminus \pi _{0}(T)$, which implies $%
f(\sigma _{BW}(T))=f(\sigma (T)\setminus \pi _{0}(T))$. \ Therefore 
\begin{eqnarray*}
f(\sigma _{BW}(T)) &=&f(\sigma (T)\setminus \pi _{0}(T)) \\
&=&\sigma (f(T))\setminus \pi _{0}(f(T))=\sigma _{BW}(f(T)).
\end{eqnarray*}

$(\Longleftarrow )$ Suppose that $f(\sigma _{BW}(T))=\sigma _{BW}(f(T))$. \
Since $T$ is isoloid, it follows from Lemma \ref{lem217} that $f(\sigma
(T)\setminus \pi _{0}(T))=\sigma (f(T))\setminus \pi _{0}(f(T))$. \ Since $%
T\in g\mathcal{W}$, we have $\sigma _{BW}(T)=\sigma (T)\setminus \pi _{0}(T)$%
. \ Therefore 
\begin{eqnarray*}
\sigma _{BW}(f(T)) &=&f(\sigma _{BW}(T)) \\
&=&f(\sigma (T)\setminus \pi _{0}(T))=\sigma (f(T))\setminus \pi _{0}(f(T)),
\end{eqnarray*}%
and hence $f(T)\in g\mathcal{W}$.
\end{proof}
\end{lemma}

As applications of Theorems \ref{thm215} and \ref{thm218} we will obtain
below several corollaries.

\begin{corollary}
\label{cor219} Suppose that $T\in \mathcal{B}(\mathcal{H})$ is reduced by
each of its eigenspaces. \ Then $f(T)\in g\mathcal{B}$ for every $f\in
H(\sigma (T))$. \ In particular, $T\in g\mathcal{B}$.
\end{corollary}

\begin{proof}
By Theorem \ref{thm215} we have $\sigma _{BW}(T)=\sigma _{BB}(T)$, so that $%
T\in g\mathcal{B}$ by Theorem \ref{thm21}. \ On the other hand, since $T$ is
reduced by each of its eigenspaces, $i(T-\lambda )i(T-\mu )\geq 0\ $for all $%
\lambda ,\mu \in \mathbb{C}\setminus \sigma _{BF}(T)$. \ It follows that $%
T\in \mathcal{P(X)}$, so Theorem \ref{thm210} implies that 
\begin{equation*}
\sigma _{BW}(f(T))=f(\sigma _{BW}(T))=f(\sigma _{BB}(T))=\sigma _{BB}(f(T)).
\end{equation*}%
Hence $f(T)\in g\mathcal{B}$.
\end{proof}

In Example \ref{ex213} we already noticed that the generalized Weyl's
theorem does not transfer to or from adjoints. \ However, we have:

\begin{corollary}
\label{cor220} Suppose that $T\in \mathcal{B}(\mathcal{H})$ is reduced by
each of its eigenspaces, and assume that $\sigma (T)$ has no isolated
points. \ Then $T,T^{\ast }\in g\mathcal{W}$. \ Moreover, if $f\in H(\sigma
(T))$ then $f(T)\in g\mathcal{W}$.
\end{corollary}

\begin{proof}
We first show that $T\in g\mathcal{W}$. \ Since $T$ is reduced by each of
its eigenspaces, it follows from Theorem \ref{thm215} that $T\in g\mathcal{B}
$. \ By Theorem \ref{thm21}, $\sigma (T)\setminus \sigma _{BW}(T)\subseteq
\pi _{0}(T)$. \ But $\operatorname*{iso}\;\sigma (T)=\emptyset $, hence $\pi
_{0}(T)=\emptyset $, which implies $\sigma _{BW}(T)=\sigma (T)$. \
Therefore, $T\in g\mathcal{W}$. \ On the other hand, observe that 
\begin{equation*}
\sigma (T^{\ast })=\overline{\sigma (T)},\ \sigma _{BW}(T^{\ast })=\overline{%
\sigma _{BW}(T)},
\end{equation*}%
and 
\begin{equation*}
\pi _{0}(T^{\ast })=\overline{\pi _{0}(T)}=\emptyset .
\end{equation*}%
Hence $T^{\ast }\in g\mathcal{W}$. \ Let $f\in H(\sigma (T))$. \ Since $T$
is reduced by each of its eigenvalues, $i(T-\lambda )i(T-\mu )\geq 0\ $\ for
all $\lambda ,\mu \in \mathbb{C}\setminus \sigma _{BF}(T)$. \ Therefore $%
\sigma _{BW}(f(T))=f(\sigma _{BW}(T))$ by Theorem \ref{thm210}. \ But $%
\sigma (T)$ has no isolated points, hence $T$ is isoloid. \ It follows from
Theorem \ref{thm218} that generalized Weyl's theorem holds for $f(T)$.
\end{proof}

For the next result, we recall that an operator $T$ is called \textit{%
reduction-isoloid} if the restriction of $T$ to every reducing subspace is
isoloid; it is well known that hyponormal operators are reduction-isoloid %
\cite{Sta}.

\begin{corollary}
\label{cor221} Suppose that $T\in \mathcal{B}(\mathcal{H})$ is both
reduction-isoloid and reduced by each of its eigenspaces. \ Then $f(T)\in g%
\mathcal{W}$ for every $f\in H(\sigma (T))$.
\end{corollary}

\begin{proof}
We first show that $T\in g\mathcal{W}$. \ In view of Theorem \ref{thm215},
it suffices to show that $\pi _{0}(T)\subseteq \sigma (T)\setminus \sigma
_{BW}(T)$. \ Suppose that $\lambda \in \pi _{0}(T)$. \ Then, with the
preceding notations, 
\begin{equation*}
\lambda \in \pi _{0}(T_{1})\cap \lbrack \operatorname*{iso}\;\sigma (T_{2})\cup
\rho (T_{2})].
\end{equation*}%
If $\lambda \in \operatorname*{iso}\;\sigma (T_{2})$, then since $T_{2}$ is
isoloid we have $\lambda \in \sigma _{p}(T_{2})$. \ But $\sigma
_{p}(T_{2})=\emptyset $, hence we must have $\lambda \in \pi _{0}(T_{1})\cap
\rho (T_{2})$. \ Since $T_{1}$ is normal, $T_{1}\in g\mathcal{W}$. \ Hence $%
T_{1}-\lambda $ is $B$-Weyl and so is $T-\lambda $, which implies $\lambda
\in \sigma (T)\setminus \sigma _{BW}(T)$. \ Therefore $\pi _{0}(T)\subseteq
\sigma (T)\setminus \sigma _{BW}(T)$, and hence $T\in g\mathcal{W}$. \ Now,
let $f\in H(\sigma (T))$. \ Since $T$ is reduced by each of its eigenspaces, 
$i(T-\lambda )i(T-\mu )\geq 0\ $for all $\lambda ,\mu \in \mathbb{C}%
\setminus \sigma _{BF}(T)$. \ It follows from Theorem \ref{thm210} that $%
f(\sigma _{BW}(T))=\sigma _{BW}(f(T))$. \ Therefore $f(T)\in g\mathcal{W}$
by Theorem \ref{thm218}.
\end{proof}

\section{Applications}

In \cite{Berkani2} and \cite{Berkani3}, the authors showed that the
generalized Weyl's theorem holds for normal operators. \ In this section we
extend this result to algebraically $M$-hyponormal operators and to
algebraically paranormal operators, using the results in Sections \ref{sect2}
and \ref{sect3}. \ We begin with the following definition.

\begin{definition}
\label{defMh}An operator $T\in \mathcal{B}(\mathcal{H})$ is said to be 
\textit{$M$-hyponormal}\ if there exists a positive real number $M$ such
that 
\begin{equation*}
M||(T-\lambda )x||\geq ||(T-\lambda )^{\ast }x||\quad \text{for all}\ x\in 
\mathcal{H}\text{, }\lambda \in \mathbb{C}.
\end{equation*}%
We say that $T\in \mathcal{B}(\mathcal{H})$ is \textit{algebraically $M$%
-hyponormal} if there exists a nonconstant complex polynomial $p$ such that $%
p(T)$ is $M$-hyponormal.
\end{definition}

The following implications hold: 
\begin{equation*}
\text{hyponormal}\Longrightarrow \text{$M$-hyponormal}\Longrightarrow \text{%
algebraically $M$-hyponormal}.
\end{equation*}%
The following result follows from Definition \ref{defMh} and some well known
facts about $M$-hyponormal operators.

\begin{lemma}
(i) If $T$ is algebraically $M$-hyponormal then so is $T-\lambda $ for every 
$\lambda \in \mathbb{C}$.\newline
(ii) If $T$ is algebraically $M$-hyponormal and $\mathcal{M}\subseteq 
\mathcal{H}$ is invariant under $T$, then $T|\mathcal{M}$ is algebraically $%
M $-hyponormal.\newline
(iii) If $T$ is $M$-hyponormal, then $N(T-\lambda )\subseteq N(T-\lambda
)^{\ast }$ for every $\lambda \in \mathbb{C}$.\newline
(iv) Every quasinilpotent $M$-hyponormal operator is the zero operator.
\end{lemma}

In \cite{Arora}, Arora and Kumar proved that Weyl's theorem holds for every $%
M$-hyponormal operator. \ We shall show that the generalized Weyl's theorem
holds for algebraically $M$-hyponormal operators. \ To do this, we need
several preliminary results.

\begin{lemma}
\label{lem32}Let $T\in \mathcal{B}(\mathcal{H})$ be $M$-hyponormal, let $%
\lambda \in \mathbb{C}$, and assume that $\sigma (T)=\{\lambda \}$. \ Then $%
T=\lambda $.

\begin{proof}
Since $T$ is $M$-hyponormal, $T-\lambda $ is also $M$-hyponormal. \ Since $%
T-\lambda $ is quasinilpotent, (iv) above implies that $T-\lambda =0$.
\end{proof}

\begin{lemma}
\label{lem33} Let $T\in \mathcal{B}(\mathcal{H})$ be a quasinilpotent
algebraically $M$-hyponormal operator. \ Then $T$ is nilpotent.
\end{lemma}

\begin{proof}
Let $p$ be a nonconstant polynomial such that $p(T)$ is $M$-hyponormal. \
Since $\sigma (p(T))=p(\sigma (T))$, the operator $p(T)-p(0)$ is
quasinilpotent. \ It follows from Lemma \ref{lem32} that $c\ T^{m}(T-\lambda
_{1})\cdots (T-\lambda _{n})\equiv p(T)-p(0)=0$. \ Since $T-\lambda _{i}$ is
invertible for every $\lambda _{i}\neq 0$, we must have $T^{m}=0$.
\end{proof}
\end{lemma}

It is well known that every $M$-hyponormal operator is isoloid. \ We can
extend this result to the algebraically $M$-hyponormal operators.

\begin{lemma}
\label{lem34} Let $T\in \mathcal{B}(\mathcal{H})$ be an algebraically $M$%
-hyponormal operator. \ Then $T$ is isoloid.

\begin{proof}
Let $\lambda $ be an isolated point of $\sigma (T)$. \ Using the spectral
projection $P:=\frac{1}{2\pi i}\int_{\partial B}(\mu -T)^{-1}d\mu $, where $%
B $ is a closed disk of center $\lambda $ which contains no other points of $%
\sigma (T)$, we can represent $T$ as the direct sum $T=T_{1}\oplus T_{2}$,$\ 
$where$\ \sigma (T_{1})=\{\lambda \}\ $and$\ \sigma (T_{2})=\sigma
(T)\setminus \{\lambda \}$. $\ $Since $T$ is algebraically $M$-hyponormal, $%
p(T)$ is $M$-hyponormal for some nonconstant polynomial $p$. \ Since $\sigma
(T_{1})=\{\lambda \}$, $\sigma (p(T_{1}))=p(\sigma (T_{1}))=\{p(\lambda )\}$%
. \ Therefore $p(T_{1})-p(\lambda )$ is quasinilpotent. \ Since $p(T_{1})$
is $M$-hyponormal, it follows from Lemma \ref{lem32} that $%
p(T_{1})-p(\lambda )=0$. \ Put $q(z):=p(z)-p(\lambda )$. \ Then $q(T_{1})=0$%
, and hence $T_{1}$ is algebraically $M$-hyponormal. \ Since $T_{1}-\lambda $
is quasinilpotent and algebraically $M$-hyponormal, it follows from Lemma %
\ref{lem33} that $T_{1}-\lambda $ is nilpotent. \ Therefore $\lambda \in
\sigma _{p}(T_{1})$, and hence $\lambda \in \sigma _{p}(T)$. \ This shows
that $T$ is isoloid.
\end{proof}

\begin{lemma}
\label{lem35} Let $T\in \mathcal{B}(\mathcal{H})$ be an algebraically $M$%
-hyponormal operator. \ Then $T$ has finite ascent. \ In particular, every
algebraically $M$-hyponormal operator has SVEP.
\end{lemma}

\begin{proof}
Suppose $p(T)$ is $M$-hyponormal for some nonconstant polynomial $p$. \
Since $M$-hyponormality is translation-invariant, we may assume $p(0)=0$. \
If $p(\lambda )\equiv a_{0}\lambda ^{m}$, then $N(T^{m})=N(T^{2m})$ because $%
M$-hyponormal operators are of ascent 1. \ Thus we write $p(\lambda )\equiv
a_{0}\,\lambda ^{m}(\lambda -\lambda _{1})\cdots (\lambda -\lambda _{n})$ ($%
m\neq 0$; $\lambda _{i}\neq 0$ for $1\leq i\leq n$). \ We then claim that 
\begin{equation}
N(T^{m})=N(T^{m+1}).  \label{eq351}
\end{equation}%
To show (\ref{eq351}), let $0\neq x\in N(T^{m+1})$. \ Then we can write 
\begin{equation*}
p(T)x=(-1)^{n}\,a_{0}\,\lambda _{1}\cdots \lambda _{n}\,T^{m}x.
\end{equation*}%
Thus we have 
\begin{align*}
|a_{0}\lambda _{1}\cdots \lambda _{n}|^{2}||T^{m}x||^{2}& =(p(T)x,\ p(T)x) \\
& \leq ||p(T)^{\ast }p(T)x||\,||x|| \\
& \leq M||p(T)^{2}x||\,||x||\quad \text{(because $p(T)$ is $M$-hyponormal)}
\\
& =M||a_{0}^{2}\,(T-\lambda _{1}I)^{2}\cdots (T-\lambda
_{n}I)^{2}T^{2m}x||\,||x|| \\
& =0,
\end{align*}%
which implies $x\in N(T^{m})$. \ Therefore $N(T^{m+1})\subseteq N(T^{m})$
and the reverse inclusion is always true. \ Since every algebraically $M$%
-hyponormal operator has finite ascent, it follows from \cite[Proposition
1.8]{Lau1} that every algebraically $M$-hyponormal operator has SVEP.
\end{proof}

\begin{theorem}
\label{thm36} Let $T\in \mathcal{B}(\mathcal{H})$ be an algebraically $M$%
-hyponormal operator. \ Then $f(T)\in g\mathcal{W}$ for every $f\in H(\sigma
(T))$.
\end{theorem}

\begin{proof}
We first show that $T\in g\mathcal{W}$. \ Suppose that $\lambda \in \sigma
(T)\setminus \sigma _{BW}(T)$. \ Then $T-\lambda $ is $B$-Weyl but not
invertible. \ Since $T$ is algebraically $M$-hyponormal, there exists a
nonconstant polynomial $p$ such that $p(T)$ is $M$-hyponormal. \ Since every
algebraically $M$-hyponormal operator has SVEP by Lemma \ref{lem35}, $T$ has
SVEP. \ It follows from Theorem \ref{thm212} that $T\in g\mathcal{B}$. \
Therefore $\sigma _{BW}(T)=\sigma _{BB}(T)$. \ But $\sigma _{BB}(T)=\sigma
_{BW}(T)\cup $ $\operatorname*{acc}\;\sigma (T)$ by Theorem \ref{thm27}, hence $%
\lambda $ is an isolated point of $\sigma (T)$. \ Since every algebraically $%
M$-hyponormal operator is isoloid by Lemma \ref{lem34}, $\lambda \in \pi
_{0}(T)$.

Conversely, suppose that $\lambda \in \pi _{0}(T)$. \ Then $\lambda $ is an
isolated eigenvalue of $T$. \ Since $\lambda $ is an isolated point of $%
\sigma (T)$, using the Riesz idempotent $E:=\frac{1}{2\pi i}\int_{\partial
D}(\mu -T)^{-1}d\mu $, where $D$ is a closed disk of center $\lambda $ which
contains no other points of $\sigma (T)$, we can represent $T$ as the direct
sum $T=T_{1}\oplus T_{2}$, where$\ \sigma (T_{1})=\{\lambda \}\ $and$\
\sigma (T_{2})=\sigma (T)\setminus \{\lambda \}$. $\ $Since $T$ is
algebraically $M$-hyponormal, $p(T)$ is $M$-hyponormal for some nonconstant
polynomial $p$. \ Since $\sigma (T_{1})=\{\lambda _{1}\}$, we have $\sigma
(p(T_{1}))=p(\sigma (T_{1}))=\{p(\lambda )\}$. \ Therefore $%
p(T_{1})-p(\lambda )$ is quasinilpotent. \ Since $p(T_{1})$ is $M$%
-hyponormal, it follows from Lemma \ref{lem32} that $p(T_{1})-p(\lambda )=0$%
. \ Define $q(z):=p(z)-p(\lambda )$. \ Then $q(T_{1})=0$, and hence $T_{1}$
is algebraically $M$-hyponormal. \ Since $T_{1}-\lambda $ is quasinilpotent
and algebraically $M$-hyponormal, it follows from Lemma \ref{lem33} that $%
T_{1}-\lambda $ is nilpotent. \ Since $T-\lambda =(T_{1}-\lambda )\oplus
(T_{2}-\lambda )$ is the direct sum of an invertible operator and a
nilpotent operator, $T-\lambda $ is $B$-Weyl. \ Hence $\lambda \in \sigma
(T)\setminus \sigma _{BW}(T)$. \ Therefore $\sigma (T)\setminus \sigma
_{BW}(T)=\pi _{0}(T)$, and hence $T\in g\mathcal{W}$.

Now let $f\in H(\sigma (T))$. \ We shall show that $\sigma
_{BW}(f(T))=f(\sigma _{BW}(T))$. \ In view of Theorem \ref{thm26} it
suffices to show that $f(\sigma _{BW}(T))\subseteq \sigma _{BW}(f(T))$. \
Suppose that $\lambda \notin \sigma _{BW}(f(T))$, and let 
\begin{equation}
f(T)-\lambda =c_{0}(T-\lambda _{1})(T-\lambda _{2})\cdots (T-\lambda
_{n})g(T),  \label{eq361}
\end{equation}%
where $c_{0},\lambda _{1},\lambda _{2},\dots ,\lambda _{n}\in \mathbb{C}$
and $g(T)$ is invertible. \ Since $\lambda \notin \sigma _{BW}(f(T))$, $%
f(T)-\lambda $ is $B$-Weyl. \ Therefore $f(T)-\lambda $ is $B$-Fredholm with
index 0. \ Since the operators on the right-hand side of (\ref{eq361})
commute, it follows from \cite[Corollary 3.3]{Berkani1} that $T-\lambda _{i}$
is $B$-Fredholm for every $i=1,2,\dots ,n$. \ Since $T$ is algebraically $M$%
-hyponormal, $T|\mathcal{M}$ is also algebraically $M$-hyponormal, where $%
\mathcal{M}$ is any closed invariant subspace of $T$. \ It follow from Lemma %
\ref{lem35} that $T$ has finite ascent. \ Hence $T|\mathcal{M}$ has also
finite ascent. \ Therefore $i(T-\lambda )\leq 0$ for every $\lambda \in 
\mathbb{C}\setminus \sigma _{BF}(T)$. \ Since $i(T-\lambda )\leq 0$ for
every $\lambda \in \mathbb{C}\setminus \sigma _{BF}(T)$ and 
\begin{equation*}
i(f(T)-\lambda )=\sum_{j=1}^{n}i(T-\lambda _{j})+i(g(T))=0,
\end{equation*}%
$T-\lambda _{i}$ is $B$-Weyl for every $i=1,2,\dots ,n$. \ Therefore $%
\lambda \notin f(\sigma _{BW}(T))$, and hence $f(\sigma _{BW}(T))\subseteq
\sigma _{BW}(f(T))$. \ Since every algebraically $M$-hyponormal operator is
isoloid by Lemma \ref{lem34}, it follows from Lemma \ref{lem217} that $%
\sigma (f(T))\setminus \pi _{0}(f(T))=f(\sigma (T)\setminus \pi _{0}(T))$. \
Hence, 
\begin{eqnarray*}
\sigma (f(T))\setminus \pi _{0}(f(T)) &=&f(\sigma (T)\setminus \pi _{0}(T))
\\
&=&f(\sigma _{BW}(T))=\sigma _{BW}(f(T)),
\end{eqnarray*}%
which implies that $f(T)\in g\mathcal{W}$.
\end{proof}

\begin{definition}
\label{defparanormal}An operator $T\in \mathcal{B}(\mathcal{H})$ is said to
be \textit{paranormal} if 
\begin{equation*}
||Tx||^{2}\leq ||T^{2}x||\quad \text{for all }x\in \mathcal{H}\text{, }%
||x||=1.
\end{equation*}%
We say that $T\in \mathcal{B}(\mathcal{H})$ is \textit{algebraically
paranormal} if there exists a nonconstant complex polynomial $p$ such that $%
p(T)$ is paranormal.
\end{definition}
\end{lemma}

The following implications hold: 
\begin{eqnarray*}
\text{hyponormal} &\Longrightarrow &\text{$p$-hyponormal} \\
&\Longrightarrow &\text{paranormal}\Longrightarrow \text{algebraically
paranormal}.
\end{eqnarray*}%
The following facts follow from Definition \ref{defparanormal} and some well
known facts about paranormal operators.

\begin{lemma}
(i) If $T\in \mathcal{B}(\mathcal{H})$ is algebraically paranormal then so
is $T-\lambda $ for every $\lambda \in \mathbb{C}$.\newline
(ii)\ If $T\in \mathcal{B}(\mathcal{H})$ is algebraically paranormal and $%
\mathcal{M}\subseteq \mathcal{H}$ is invariant under $T$, then $T|\mathcal{M}
$ is algebraically paranormal.
\end{lemma}

In \cite{Curto2} we showed that if $T$ is an algebraically paranormal
operator then $f(T)\in \mathcal{W}$ for every $f\in H(\sigma (T))$. \ We can
now extend this result to the generalized Weyl's theorem. \ To prove this we
need several lemmas.

\begin{lemma}
\label{lem38}Let $T\in \mathcal{B}(\mathcal{H})$ be $B$-Fredholm. \ The
following statements are equivalent:\newline
(i) $T$ does not have SVEP at $0$;\newline
(ii) $a(T)=\infty $;\newline
(iii) $0\in $ $\operatorname*{acc}\;\sigma _{p}(T)$.

\begin{proof}
Suppose that $T$ is $B$-Fredholm. \ It follows from \cite[Theorem 2.7]%
{Berkani1} that $T$ can be decomposed as 
\begin{equation*}
T=T_{1}\oplus T_{2}\,\;\text{(}T_{1}\text{ Fredholm, }T_{2}\text{ nilpotent)}%
.
\end{equation*}%
(i)$\Longleftrightarrow $(ii): Suppose that $T$ does not have SVEP at $0$. \
Since $T_{2}$ is nilpotent, $T_{2}$ has SVEP. \ Therefore $T_{1}$ does not
have SVEP. \ Since $T_{1}$ is Fredholm, it follows from \cite[Theorem 2.6]%
{Aie} that $a(T)=\infty $.

Conversely, suppose that $a(T)=\infty $. \ Since $T_{2}$ is nilpotent, $%
T_{2} $ has finite ascent. \ Therefore $a(T_{1})=\infty $. \ But $T_{1}$ is
Fredholm, hence $T_{1}$ does not have SVEP by \cite[Theorem 2.6]{Aie}.

(i)$\Longleftrightarrow $(iii): Suppose that $T$ does not have SVEP at $0$.
\ Then $T_{1}$ does not have SVEP. \ Since $T_{1}$ is Fredholm, it follows
from \cite[Theorem 2.6]{Aie} that $0\in \operatorname*{acc}\;\sigma _{p}(T_{1})$.
\ Therefore $0\in \operatorname*{acc}\;\sigma _{p}(T)$.

Conversely, suppose that $0\in \operatorname*{acc}\;\sigma _{p}(T)$. \ Since $%
T_{2} $ is nilpotent, $0\in \operatorname*{acc}\;\sigma _{p}(T_{1})$. \ But $T_{1}$
is Fredholm, hence $T_{1}$ does not have SVEP by \cite[Theorem 2.6]{Aie}. \
Therefore $T$ does not have SVEP.
\end{proof}

\begin{corollary}
Suppose that $T\in \mathcal{B}(\mathcal{H})$ is $B$-Fredholm with $i(T)>0$.
\ Then $T$ does not have SVEP at $0$.
\end{corollary}

\begin{proof}
Suppose that $T$ is $B$-Fredholm with $i(T)>0$. \ Then by \cite[Theorem 2.7]%
{Berkani1}, $T$ can be decomposed by 
\begin{equation*}
T=T_{1}\oplus T_{2}\,\;\text{(}T_{1}\text{ Fredholm, }T_{2}\text{ nilpotent)}%
.
\end{equation*}%
Moreover, $i(T)=i(T_{1})$. \ But $i(T)>0$, hence $i(T_{1})>0$. \ Since $%
T_{1} $ is Fredholm, it follows from \cite[Corollary 11]{Fin} that $T_{1}$
does not have SVEP at $0$. \ Therefore $T$ does not have SVEP at $0$.
\end{proof}

\begin{theorem}
Suppose that $T\in \mathcal{B}(\mathcal{H})$ is $B$-Fredholm. \ Then 
\begin{equation*}
T^{\ast }\text{ does not have SVEP at }0\ \Longleftrightarrow \ d(T)=\infty .
\end{equation*}%
Moreover, if $T$ and $T^{\ast }$ have SVEP at $0$ then $T$ is $B$-Fredholm
with index $0$.
\end{theorem}

\begin{proof}
Since $T$ is $B$-Fredholm, $T$ can be decomposed by 
\begin{equation*}
T=T_{1}\oplus T_{2}\,\;\text{(}T_{1}\text{ Fredholm, }T_{2}\text{ nilpotent)}%
.
\end{equation*}%
But $T_{1}$ is Fredholm if and only if $T_{1}^{\ast }$ is Fredholm, hence $T$
is $B$-Fredholm if and only if $T^{\ast }$ is $B$-Fredholm. \ Since $T_{1}$
is Fredholm, $a(T_{1})=d(T_{1}^{\ast })$. \ Also, since $T_{2}$ is
nilpotent, $a(T_{2})=d(T_{2})=a(T_{2}^{\ast })=d(T_{2}^{\ast })$. \ It
follows from \cite[Theorem 6.1]{A.E.Taylor} that 
\begin{equation*}
\begin{split}
a(T^{\ast })& =a(T_{1}^{\ast }\oplus T_{2}^{\ast }) \\
& =\max \{a(T_{1}^{\ast }),a(T_{2}^{\ast })\} \\
& =\max \{d(T_{1}),d(T_{2})\} \\
& =d(T_{1}\oplus T_{2}) \\
& =d(T).
\end{split}%
\end{equation*}%
Therefore by Lemma \ref{lem38}, 
\begin{equation*}
T^{\ast }\ \text{does not have SVEP at }0\Longleftrightarrow \ a(T^{\ast
})=\infty \Longleftrightarrow \ d(T)=\infty .
\end{equation*}%
Moreover, suppose that $T$ and $T^{\ast }$ have SVEP at $0$. \ Then by Lemma %
\ref{lem38}, $a(T)=d(T)<\infty $, and hence $T$ is $B$-Fredholm with index $%
0 $.
\end{proof}

\begin{lemma}
\label{lem311}\cite[Lemmas 2.1, 2.2, 2.3]{Curto2}) Let $T\in \mathcal{B(H)}$
be an algebraically paranormal operator. Then\newline
(i) If $\sigma (T)=\{\lambda \}$, then $T=\lambda $;\newline
(ii) If $T$ is quasinilpotent, then it is nilpotent;\newline
(iii) $T$ is isoloid.
\end{lemma}

\begin{theorem}
\label{algpara}Let $T\in \mathcal{B(H)}$ be an algebraically paranormal
operator. \ Then $f(T)\in g\mathcal{W}$ $\ $for every $f\in H(\sigma (T))$.
\end{theorem}

\begin{proof}
We first show that $T\in g\mathcal{W}$. \ Suppose that $\lambda \in \sigma
(T)\setminus \sigma _{BW}(T)$. \ Then $T-\lambda $ is $B$-Weyl but not
invertible. \ Since $T$ is an algebraically paranormal operator, there
exists a nonconstant polynomial $p$ such that $p(T)$ is paranormal. \ Since
every paranormal operator has SVEP, $p(T)$ has SVEP. \ Therefore $T$ has
SVEP. \ It follows from Theorem \ref{thm212} that $T\in g\mathcal{B}$. \
Therefore $\sigma _{BW}(T)=\sigma _{BB}(T)$. \ But $\sigma _{BB}(T)=\sigma
_{BW}(T)\cup \operatorname*{acc}\;\sigma (T)$ by Theorem \ref{thm27}, hence $%
\lambda $ is an isolated point of $\sigma (T)$. \ Since every algebraically
paranormal operator is isoloid by Lemma \ref{lem311}, $\lambda \in \pi
_{0}(T)$.

Conversely, suppose that $\lambda \in \pi _{0}(T)$. \ Let $P:=\frac{1}{2\pi i%
}\int_{\partial D}(\mu -T)^{-1}d\mu $ be the associated Riesz idempotent,
where $D$ is an open disk of center $\lambda $ which contains no other
points of $\sigma (T)$, we can represent $T$ as the direct sum $%
T=T_{1}\oplus T_{2}$, where$\ \sigma (T_{1})=\{\lambda \}\ $and$\ \sigma
(T_{2})=\sigma (T)\setminus \{\lambda \}$. $\ $Now we consider two cases:

\textbf{Case I.} \ Suppose that $\lambda =0$. \ Then $T_{1}$ is
algebraically paranormal and quasinilpotent. \ It follows from Lemma \ref%
{lem311} that $T_{1}$ is nilpotent. \ Therefore $T$ is the direct sum of an
invertible operator and nilpotent, and hence $T$ is $B$-Weyl by \cite[Lemma
4.1]{Berkani2}. \ Thus, $0\in \sigma (T)\setminus \sigma _{BW}(T)$.

\textbf{Case II. \ }Suppose that $\lambda \neq 0$. \ Since $T$ is
algebraically paranormal, $p(T)$ is paranormal for some nonconstant
polynomial $p$. \ Since $\sigma (T_{1})=\{\lambda _{1}\}$, we have $\sigma
(p(T_{1}))=p(\sigma (T_{1}))=\{p(\lambda )\}$. \ Therefore $%
p(T_{1})-p(\lambda )$ is quasinilpotent. \ Since $p(T_{1})$ is paranormal,
it follows from Lemma \ref{lem311} that $p(T_{1})-p(\lambda )=0$. \ Define $%
q(z):=p(z)-p(\lambda )$. \ Then $q(T_{1})=0$, and hence $T_{1}$ is
algebraically paranormal. \ Since $T_{1}-\lambda $ is quasinilpotent and
algebraically paranormal, it follows from Lemma \ref{lem311} that $%
T_{1}-\lambda $ is nilpotent. \ Since $T-\lambda =%
\begin{pmatrix}
T_{1}-\lambda & 0 \\ 
0 & T_{2}-\lambda%
\end{pmatrix}%
$ is the direct sum of an invertible operator and nilpotent, $T-\lambda $ is 
$B$-Weyl. \ Therefore $\lambda \in \sigma (T)\setminus \sigma _{BW}(T)$. \
Thus $T\in g\mathcal{W}$.

Now we claim that $\sigma _{BW}(f(T))=f(\sigma _{BW}(T))$ for every $f\in
H(\sigma (T))$. \ Let $f\in H(\sigma (T))$. \ Since $\sigma
_{BW}(f(T))\subseteq f(\sigma _{BW}(T))$ with no other restriction on $T$ by
Theorem \ref{thm26}, it suffices to show that $f(\sigma _{BW}(T))\subseteq
\sigma _{BW}(f(T))$. \ Suppose that $\lambda \notin \sigma _{BW}(f(T))$. \
Then $f(T)-\lambda $ is $B$-Weyl and 
\begin{equation}
f(T)-\lambda \equiv c_{0}(T-\lambda _{1})(T-\lambda _{2})\cdots (T-\lambda
_{n})g(T),  \label{eq3121}
\end{equation}%
where $c_{0},\lambda _{1},\lambda _{2},\dots ,\lambda _{n}\in \mathbb{C}$
and $g(T)$ is invertible. \ Since the operators on the right-hand side of (%
\ref{eq3121}) commute, every $T-\lambda _{i}$ is $B$-Fredholm by %
\cite[Corollary 3.3]{Berkani1}. \ Since $T$ is algebraically paranormal, $T$
has SVEP. \ It follows from Lemma \ref{lem38} that $i(T-\lambda _{i})\leq 0$
(all $i=1,2,\dots ,n$). \ Therefore $\lambda \notin f(\sigma _{BW}(T))$, and
hence $\sigma _{BW}(f(T))=f(\sigma _{BW}(T))$. \ Since $T$ is algebraically
paranormal, it follows from Lemma \ref{lem311} that $T$ is isoloid. \
Therefore by Lemma \ref{lem217}, 
\begin{equation*}
\sigma (f(T))\setminus \pi _{0}(f(T))=f(\sigma (T)\setminus \pi _{0}(T)).
\end{equation*}%
Hence 
\begin{equation*}
\sigma (f(T))\setminus \pi _{0}(f(T))=f(\sigma (T)\setminus \pi
_{0}(T))=f(\sigma _{BW}(T))=\sigma _{BW}(f(T)),
\end{equation*}%
which implies that $f(T)\in g\mathcal{W}$.
\end{proof}
\end{lemma}


\begin{thebibliography}{99}
\bibitem{Aie} P. Aiena and O. Monsalve, Operators which do not have the
single valued extension property\textit{,} \textit{J. Math. Anal. Appl.} 
\textbf{250} (2000), 435--448.

\bibitem{Arora} S.C. Arora and R. Kumar, $M$-hyponormal operators\textit{,} 
\textit{Yokohama Math. J}. \textbf{28} (1980), 41--44.

\bibitem{Ber1} S.K. Berberian, An extension of Weyl's theorem to a class of
not necessarily normal operators\textit{,} \textit{Michigan Math. J.} 
\textbf{16} (1969), 273--279.

\bibitem{Ber2} S.K. Berberian, The Weyl spectrum of an operator,\textit{\
Indiana Univ. Math. J}. \textbf{20} (1970), 529--544.

\bibitem{Berkani1} M. Berkani, On a class of quasi-Fredholm operators, 
\textit{Integral Equations Operator Theory,} \textbf{34} (1999), 244--249.

\bibitem{Berkani2} M. Berkani, Index of $B$-Fredholm operators and
generalization of a Weyl theorem, \textit{Proc. Amer. Math. Soc.} \textbf{130%
} (2002), 1717--1723.

\bibitem{Berkani3} M. Berkani, $B$-Weyl spectrum and poles of the resolvent, 
\textit{J. Math. Anal. Appl.} \textbf{272} (2002), 596--603.

\bibitem{Berkani4} M. Berkani and J.J. Koliha, Weyl type theorems for
bounded linear operators, \textit{Acta Sci. Math. (Szeged) }\textbf{69}
(2003), 359--376.

\bibitem{Berkani6} M. Berkani and M. Sarih, On semi $B$-Fredholm operators, 
\textit{Glasgow Math. J.} \textbf{43} (2001), 457--465.

\bibitem{Co} L.A. Coburn, Weyl's theorem for nonnormal operators, \textit{%
Michigan Math. J.} \textbf{13} (1966), 285--288.

\bibitem{Curto1} R.E. Curto and Y.M. Han, Weyl's theorem, $a$-Weyl's
theorem, and local spectral theory, \textit{J. London Math. Soc.} (2) 
\textbf{67} (2003), 499--509.

\bibitem{Curto2} R.E. Curto and Y.M. Han, Weyl's theorem holds for
algebraically paranormal operators, \textit{Integral Equations Operator
Theory} \textbf{47} (2003), 307--314.

\bibitem{Fin} J.K. Finch, The single valued extension property on a Banach
space, \textit{Pacific J. Math.} \textbf{58} (1975), 61--69.

\bibitem{Har1} R.E. Harte, Fredholm, Weyl and Browder theory, \textit{Proc.
Royal Irish Acad. }\textbf{85A} (1985), 151--176.

\bibitem{Har2} R.E. Harte, \textit{Invertibility and Singularity for Bounded
Linear Operators},\textit{\ }Marcel Dekker, New York, 1988.

\bibitem{Har3} R.E. Harte and W.Y. Lee, Another note on Weyl's theorem, 
\textit{Trans. Amer. Math. Soc.} \textbf{349} (1997), 2115--2124.

\bibitem{Lau1} K.B. Laursen, Operators with finite ascent, \textit{Pacific
J. Math. }\textbf{152} (1992), 323--336.

\bibitem{Lau2} K.B. Laursen and M.M. Neumann, \textit{An Introduction to
Local Spectral Theory}, London Mathematical Society Monographs New Series
20, Clarendon Press, Oxford, 2000.

\bibitem{Oberai2} K.K. Oberai, On the Weyl spectrum (II), \textit{Illinois
J. Math. }\textbf{21} (1977), 84--90.

\bibitem{Sta} J. Stampfli, Hyponormal operators, \textit{Pacific J. Math.}
12(1962), 1453-1458.

\bibitem{A.E.Taylor} A.E. Taylor, Theorems on ascent, descent, nullity and
defect of linear operators, \textit{Math. Ann. }\textbf{163} (1966), 18--49.

\bibitem{Weyl} H. Weyl, \"{U}ber beschr\"{a}nkte quadratische Formen, deren
Differenz vollsteig ist, \textit{Rend. Circ. Mat. Palermo }\textbf{27}
(1909), 373--392.
\end{thebibliography}
\end{document}